\documentclass[a4paper,11pt]{article}
\usepackage {amsthm, amssymb, amsmath, amsfonts, enumerate,epsfig,comment}

\newtheorem{de}{Definition}
\newtheorem{te}[de]{Theorem}
\newtheorem{pr}[de]{Proposition}
\newtheorem{lm}[de]{Lemma}
\newtheorem{co}[de]{Corollary}

\newtheorem{cj}[de]{Conjecture}

\newcommand{\N}{{\mathbb N}}
\newcommand{\R}{{\mathbb R}}
\newcommand{\Q}{{\mathbb Q}}
\newcommand{\C}{{\mathbb C}}
\newcommand{\Z}{{\mathbb Z}}
\newcommand{\F}{{\mathbb F}}
\DeclareMathOperator{\ord}{ord}
\DeclareMathOperator{\NV}{NV}
\DeclareMathOperator{\ac}{ac}
\DeclareMathOperator{\real}{Re}
\DeclareMathOperator{\mult}{mult}
\DeclareMathOperator{\proj}{pr}
\DeclareMathOperator{\id}{id}

\begin{document}
\title{The holomorphy conjecture for \\ nondegenerate surface singularities}
\author{Wouter Castryck, Denis Ibadula and Ann Lemahieu \footnote{
 The
research was partially supported by MCI-Spain grant MTM2010-21740-C02, by the ANR `SUSI'
project (ANR-12-JS01-0002-01) and by the research project G093913N of the Research Foundation - Flanders (FWO).}
\date{}}
\maketitle {\footnotesize \emph{\noindent \textbf{Abstract.---} The
holomorphy conjecture states roughly that Igusa's zeta function
associated to a hypersurface and a character is holomorphic on
$\mathbb{C}$ whenever the order of the character does not divide the
order of any eigenvalue of the local monodromy of the hypersurface.
In this article we prove the holomorphy conjecture for surface singularities which are nondegenerate over $\C$ w.r.t.\ their Newton polyhedron.
In order to provide relevant eigenvalues of monodromy, we first show a relation between the normalized volume (which appears in the formula of Varchenko for the zeta function of monodromy) of faces in a simplex in arbitrary dimension.
We then study some specific character sums that show up when dealing with false poles. In contrast with the context of the trivial character, we here need to show fakeness of certain poles in addition to the candidate poles contributed by $B_1$-facets.
}}
\section{Introduction}
Let $K$ be a finite extension of the field of $p$-adic numbers $\mathbb{Q}_p$. Let $R$ be the valuation ring of $K$ and $P$ its maximal ideal. Suppose that the residue field $R/P$ has cardinality $q$. For $z\in K$, let $\ord(z)\in \mathbb{Z}\cup\left\{\infty\right\}$ denote its valuation, $|z|=q^{-\ord(z)}$ its absolute value and $\ac(z)=z\pi^{-\ord(z)}$ its angular component, where $\pi$ is a fixed uniformizing parameter for $R$.

Let $f(\underline{x})$, $\underline{x}:=\left(x_1,\ldots,x_n\right)$, be a non-constant polynomial over $K$, and $\chi: R^{\times}\rightarrow \mathbb{C}^{\times}$ a multiplicative character of $R^{\times}$, i.e.\ a homomorphism with finite image. We formally put $\chi(0)=0$.
Let $Z_{f,0}(\chi,K,s)$, resp.\ $Z_{f}(\chi,K,s)$, be the corresponding \emph{local Igusa zeta function}, resp.\ \emph{global Igusa zeta function}, i.e.\ the meromorphic continuation to $\mathbb{C}$ of the integral function
\[Z_0(s)=\int_{P^n}\chi\left(\ac(f(\underline{x}))\right)|f(\underline{x})|^s|d(\underline{x})|,\quad \mbox{resp.\ } Z(s)=\int_{R^n}\chi\left(\ac(f(\underline{x}))\right)|f(\underline{x})|^s|d(\underline{x})|,\] for $s\in \mathbb{C}$ with $\real(s)>0$, where $|d(\underline{x})|=|dx_1\wedge\ldots \wedge dx_n|$ denotes the Haar measure on $K^n$ normalized such that the measure of $R^n$ is $1$.
%Analogously, let $Z_{f}(\chi,K,s)$ be \emph{the global Igusa zeta function}, i.e. the meromorphic continuation to $\mathbb{C}$ of the integral function
%\[Z(s)=\int_{R^n}\chi\left(ac(f(\underline{x}))\right)|f(\underline{x})|^s|d(\underline{x})|, \qquad Re (s)>0.\]

For $f$ a polynomial over $R$, the local and global Igusa zeta function can be described in terms of solutions of congruences. For $i \in \N_{>0}$ and $u \in R/P^i$, let $M_{0,i}(u)$ and $M_i(u)$ be the numbers of solutions of $f(\underline{x}) \equiv u \mbox{ mod } P^i$ in $(P/P^i)^n$ resp.\ $(R/P^i)^n$. Let $c$ be the conductor of $\chi$, i.e.\ the smallest $a \in \N_{> 0}$ such that $\chi$ is trivial on $1+P^a$. Then
$$Z_0(s) =\sum_{i=0}^{\infty}\sum_{u \in (R/P^c)^{\times}}\chi(u)M_{0,i+c}(\pi^i u)q^{-n(i+c)}q^{-is}, \qquad \mbox{and}$$
$$Z(s) =\sum_{i=0}^{\infty}\sum_{u \in (R/P^c)^{\times}}\chi(u)M_{i+c}(\pi^i u)q^{-n(i+c)}q^{-is}.$$
Igusa showed that these functions are rational functions in $q^{-s}$ and he gave a formula for $Z_{f,0}(\chi,K,s)$ and $Z_{f}(\chi,K,s)$ in terms of an embedded resolution $(Y,h)$ of $f^{-1}\{0\}$ over $K$ (see \cite{I75}).
Let $E_j$, $j\in T$, be the (reduced) irreducible components of $h^{-1}(f^{-1}\left\{0\right\})$, and let $N_j$, resp.\ $\nu_j-1$, be the multiplicity of $E_j$ in the divisor of $f\circ h$, resp.\ $h^*(dx_1 \wedge \ldots \wedge dx_n)$ on $Y$. Then the poles of $Z_{f,0}(\chi,K,s)$ and $Z_{f}(\chi,K,s)$ are among the values
\begin{eqnarray} \label{eqpole}
s=\frac{-\nu_j}{N_j}+\frac{2k\pi i}{N_j \log(q)}, \quad k \in \Z, j \in T,
\end{eqnarray}
for which the order of $\chi$ divides $N_j$.

Let now $f\in F[\underline{x}]$, with $F \subset \C$ a number field, and let $K$ be a non-archimedean completion of $F$, i.e.\ a completion w.r.t.\ a finite prime. Let $R$ be its valuation ring and let $\chi : R^\times \rightarrow \C^\times$ be a character.
Then the poles of $Z_{f,0}(\chi,K,s)$ and $Z_f(\chi,K,s)$ seem to be related to various invariants in singularity theory, such as the eigenvalues of monodromy and the roots of the Bernstein-Sato polynomial (see for example \cite{D91-2}) and such as the jumping numbers (see for example \cite{ST}).
In this article we explore another connection conjectured by Denef, called the holomorphy conjecture.
It follows from (\ref{eqpole}) that when the order of $\chi$ divides no $N_j$ at all, then the zeta functions $Z_{f,0}(\chi,K,s)$ and $Z_{f}(\chi,K,s)$ are holomorphic on $\mathbb{C}$.
Now, the $N_j$ are not intrinsically associated to $f^{-1}\left\{0\right\}$; however the order (as root of unity) of any eigenvalue of the local monodromy on $f^{-1}\left\{0\right\}$ divides some $N_j$, and those eigenvalues are intrinsic invariants of $f^{-1}\left\{0\right\}$. This observation inspired Denef to propose the following (\cite[Conjecture 4.4.2]{D91-2}):

\begin{cj}[\textbf{Holomorphy conjecture}] For almost all non-archimedean completions $K$ of $F$ (i.e.\ for all except a finite number)
and all characters $\chi$, the local (resp.\ global) Igusa zeta function $Z_{f,0}(\chi,K,s)$ (resp.\ $Z_{f}(\chi,K,s)$) is holomorphic, unless the order of $\chi$ divides the order of some eigenvalue of the local monodromy of $f$ at some complex point of $f^{-1}\left\{0\right\}$. \end{cj}

This conjecture has been proven by Veys in \cite[Theorem 3.1]{Veys holomorphy} for plane
curves, and in \cite{DV95} Denef and Veys got a Thom-Sebastiani type result. In \cite{holrodrveys} Rodrigues and Veys make several progresses on the holomorphy conjecture for homogeneous polynomials. Veys and the third author confirmed the conjecture for
surfaces that are general for a toric idealistic cluster (see
\cite[Theorem 24]{monodconjtoric}). In \cite{LVPhol} the holomorphy conjecture has been introduced for ideals and was proven for ideals in dimension two.

In this article we prove the holomorphy conjecture for surface singularities that are nondegenerate over $\C$ w.r.t.\ their
Newton polyhedron at the origin. In Section 2 we recall this notion, along with explicit formulas for
the zeta functions in this context.
By a formula of Varchenko the normalized volume of a face gets a key role in the search for eigenvalues of monodromy for nondegenerate singularities. In Section 3 we prove some properties on the normalized volume of faces in a simplex of arbitrary dimension. These properties might be of independent interest.
We can use them in Subsection \ref{subseigenvalues} to obtain a set of eigenvalues that is relevant for the holomorphy conjecture. Furthermore, we prove that some candidate poles of $Z_{f,0}(\chi,K,s)$ (resp.\ $Z_{f}(\chi,K,s)$) are no actual poles. In \cite{BV15} some configurations of $B_1$-facets that give rise to false poles have been treated in the context of non-trivial characters. Actually, it seems that almost all configurations of $B_1$-facets give rise to fake poles (see Subsection \ref{subsB1} for the exact statement). We even find a configuration without $B_1$-facets where we need to show that the candidate pole is a false pole. These computations rely on the study of some specific character sums (see Section \ref{sectioncharsums}). We can then complete our proof using a nondegeneracy argument (see Lemma \ref{nondegeneratestable}) which was used for the first time in \cite{LVP11}.

As a preliminary remark, we note that for the purpose of proving the holomorphy conjecture one can assume that
\begin{itemize}
\item $f$ has coefficients in
the ring of integers $\mathcal{O}_F$ of $F$; indeed, multiplying
$f$ by a constant $a \in F$ affects $Z_{f,0}(\chi,K,s)$ and $Z_{f}(\chi,K,s)$ only for the completions $K$ in which $\ord(a) \neq 0$, of which
there are finitely many,
\item $\chi$ is a non-trivial character with conductor equal to $1$; indeed,
Denef proved that for almost all non-archimedean completions $K$ of $F$, if $\chi: R^{\times}\rightarrow \mathbb{C}^{\times}$ is a multiplicative character which is non-trivial on $1+P$, then $Z_{f,0}(\chi,K,s)$ and
$Z_{f}(\chi,K,s)$ are constant on $\C$ (see \cite[Theorem 3.3]{D91-2}).
\end{itemize}

From now on we will just write $Z_{f,0}(\chi,s)$ (resp.\ $Z_{f}(\chi,s)$) for $Z_{f,0}(\chi,K,s)$ (resp.\ $Z_{f}(\chi,K,s)$).

\section {Nondegenerate singularities and their zeta functions}

\subsection{Nondegenerate singularities} \label{subsection_nondegsing}

Assume that $f(\underline{x}) \in \mathcal{O}_F[\underline{x}]$
is a non-constant polynomial satisfying $f(\underline{0})=0$. Write
\[ f(\underline{x}) = \sum_{\underline{k} \in \Z_{\geq 0}^n} c_{\underline{k}} \underline{x}^{\underline{k}}, \]
where $\underline{k}=(k_1, \ldots, k_n)$
and $\underline{x}^{\underline{k}} = x_1^{k_1} \cdot \ldots \cdot x_n^{k_n}.$ The \emph{support of $f$} is $\textnormal{supp}\, f = \{\underline{k} \in \Z_{\geq
0}^n \, | \, c_{\underline{k}} \neq 0 \}.$ The \emph{Newton polyhedron
$\Gamma_0$ of $f$ at the origin} is the convex hull in
$\mathbb{R}_{\geq 0}^n$ of
$$\bigcup_{\underline{k} \in \textnormal{{supp}}\, f} \underline{k} +
\mathbb{R}_{\geq 0}^n.$$
A \emph{facet} of the Newton polyhedron is a face
of dimension $n-1.$ For a face $\tau$ of $\Gamma_0$, one defines
the polynomial $f_\tau(\underline{x}) := \sum_{\underline{k} \in \Z^n \cap \tau} c_{\underline{k}} \underline{x}^{\underline{k}}$.

We say that the
polynomial $f$ is \emph{nondegenerate over $\C$ w.r.t.\ the compact faces of $\Gamma_0$} (resp.\ \emph{nondegenerate over $\C$ w.r.t.\ the faces of $\Gamma_0$}), if for every compact face $\tau$ (resp.\ for every face $\tau$) of $\Gamma_0$, the zero locus of $f_\tau$ has no singularities in
 $(\C^{\times})^n.$ For a fixed Newton polyhedron $\Gamma$, almost all polynomials
having $\Gamma$ as their Newton polyhedron are nondegenerate w.r.t.\ the faces of $\Gamma$ (see \cite[p.157]{AVG86}).

Let $K$ be a non-archimedean completion of $F$ with valuation ring $R$ and maximal ideal $P$,
whose residue field we denote by $\mathbb{F}_q$. Note that $\mathcal{O}_F \subset R$, so it
makes sense to consider $\bar{f}$, the polynomial over
$\mathbb{F}_q$ obtained from $f$ by reducing each of its coefficients modulo $P$.
We say that $\bar{f}$ is \emph{nondegenerate over $\F_q$ w.r.t.\ the compact faces of $\Gamma_0$} (resp.\ \emph{nondegenerate over $\F_q$ w.r.t.\ the faces of $\Gamma_0$}) if for every compact face $\tau$ (resp.\ for every face $\tau$) of $\Gamma_0$, the zero locus of $\bar{f_\tau}$ has no singularities in
 $(\F_q^\times)^n$. If $f$ is nondegenerate over $\C$
w.r.t.\ the compact faces (resp.\ the faces) of its Newton polyhedron $\Gamma_0$, then recall that $\bar{f}$
is nondegenerate over $\mathbb{F}_q$ w.r.t.\ the compact faces (resp.\ the faces) of $\Gamma_0$ for
almost all choices of $K$. Thus in order to prove the holomorphy conjecture for polynomials that are nondegenerate over $\C$,
it suffices to restrict to completions $K$ for
which moreover $\bar{f}$ is nondegenerate over the residue field $\mathbb{F}_q$.

Further on we will use the following property of nondegeneracy (\cite[Lemma 9]{LVP11}):
\begin{lm} \label{nondegeneratestable}
If a complex polynomial $f(x,y,z)$ is nondegenerate w.r.t.\
the compact faces of its Newton polyhedron at the origin, then for almost all $k \in
\C$ the polynomial $f(x,y,z-k)$ is nondegenerate w.r.t.\ the compact faces of
its Newton polyhedron at the origin. (Analogously for the
variables $x$ and $y$.)
\end{lm}

\subsection{Some combinatorial data associated to the Newton polyhedron}

Let $\Gamma_0$ be as above. For $\underline{a}= (a_1,\ldots , a_n) \in \R_{\geq 0}^n$ we put
\[ N(\underline{a}) :=
\inf_{\underline{x} \in \Gamma_0} \underline{a} \cdot \underline{x}, \qquad \nu (\underline{a}) :=
\sum_{i=1}^n a_i, \qquad F(\underline{a}) := \{\underline{x} \in \Gamma_0 \, | \, \underline{a} \cdot \underline{x}
= N(\underline{a})\}.\] All $F(\underline{a}), \underline{a} \neq \underline{0},$ are faces of $\Gamma_0.$ To a
face $\tau$ of $\Gamma_0$ we associate its dual cone $\Delta_\tau = \{\underline{a} \in
\R_{\geq 0}^n \, | \, F(\underline{a}) = \tau \}.$ It is a rational polyhedral
cone of dimension $n - \dim
\tau$.
%Note that the map $\tau \mapsto \Delta_\tau$ is inclusion-reversing.
In particular if $\tau$ is a facet
then $\Delta_\tau$ is a ray, say $\Delta_\tau = \underline{a} \R_{> 0}$
for some non-zero $\underline{a} \in \Z_{\geq 0}^n,$ and then the equation of the
hyperplane through $\tau$ is $\underline{a} \cdot \underline{x}= N(\underline{a}).$ If we demand that
$\underline{a}$ is primitive, i.e.\ that
$\textnormal{gcd}(a_1, \ldots , a_n)=1$, then this $\underline{a}$ is uniquely
defined. For a facet $\tau$ we also use the notation $N(\tau)$ - called the lattice distance of $\tau$ - and
$\nu(\tau)$, meaning respectively $N(\underline{a})$ and $\nu(\underline{a})$ for this
associated $\underline{a} \in \Z_{\geq 0}^n.$ For a general proper face $\tau$ the dual
cone $\Delta_\tau$ is strictly positively spanned by the dual cones of the facets containing $\tau$.

For a set of linearly independent vectors $\underline{a_1} , \ldots , \underline{a_r} \in \Z^n$
we define the \emph{multiplicity} $\mult(\underline{a_1}, \dots, \underline{a_r})$ as the
index of the lattice $\Z \underline{a_1} +
\ldots + \Z \underline{a_r}$ in the group of the points with integral
coordinates of the subspace of $\R^n$ generated by $\underline{a_1} , \ldots, \underline{a_r}.$
Alternatively, $\mult(\underline{a_1}, \dots, \underline{a_r})$ is equal to the greatest common divisor
of the determinants of the $(r \times r)$-matrices obtained by
omitting columns from the matrix with rows $\underline{a_1}, \ldots, \underline{a_r}.$
If $\Delta_\tau$ is a simplicial cone then by $\mult(\Delta_\tau)$ we
mean the multiplicity of its set of primitive generators. For a simplicial face $\tau$ we
write $\mult(\tau)$ for the multiplicity of its set of vertices.

\subsection{The Igusa zeta function with character for nondegenerate singularities} \label{subsigchar}

In the case where $f \in R[\underline{x}]$ is nondegenerate over $\mathbb{F}_q$ w.r.t.\ the compact faces (resp.\ the faces)
of its Newton polyhedron at the origin $\Gamma_0$, Hoornaert gave a formula (\cite[Theorem 3.4]{H02}) for the local (resp.\ global) Igusa zeta function associated
 to $f$ and $\chi$ in terms of $\Gamma_0$, which we recall.
 Hoornaert states the formula for $R = \mathbb{Z}_p$ only, but her proof generalizes word by word to our more general setting.

Recall that we assume $\chi : R^\times \rightarrow \mathbb{C}^\times$ to be non-trivial of conductor $1$.
 Let $\proj: R^{\times} \rightarrow \F_q^{\times} \cong R^{\times}/(1+P)$ be
the natural surjective homomorphism. As $\chi$ is trivial on $1+P$, there exists a unique homomorphism
$\bar{\chi}: \F_q^{\times} \rightarrow \C^{\times}$
such that $\chi=\bar{\chi} \circ \pi$.
One formally puts $\bar{\chi}(0) = 0$. Note that the order of $\chi$ divides the order of $\bar{\chi}$.
Let $f$ be a non-zero polynomial over $R$ satisfying $f(\underline{0})=0$ and let $\bar{f}$ be nondegenerate over $\F_q$ w.r.t.\ all the compact faces (resp.\ all the faces) of its Newton polyhedron $\Gamma_0$.
Let
$$L_\tau := q^{-n} \sum_{\underline{x} \in (\F_q^\times)^n} \bar{\chi}(\bar{f_\tau}(\underline{x})) \qquad \mbox{and} \qquad S(\Delta_\tau)(s):=\sum_{\underline{a} \in \Z^n \cap \Delta_\tau}q^{-\nu(\underline{a})-N(\underline{a})s}.$$
Then Hoornaert proved that the local, resp.\ global, Igusa zeta function associated to $f$ and
the non-trivial character $\chi$ can be computed as
\begin{eqnarray*}
Z_{f,0}(\chi,s) = \sum_{\substack{\tau \textnormal{ compact}\\
\textnormal{face of }\Gamma_0}} L_\tau S(\Delta_\tau)(s), \quad \mbox{ resp.\ } \quad Z_{f}(\chi,s) = \sum_{\substack{\tau \\
\textnormal{face of }\Gamma_0}} L_\tau S(\Delta_\tau)(s).
\end{eqnarray*}
%Analogously, for the global Igusa zeta function associated to $f$ and
%the non-trivial character $\chi$, Hoornaert proves that
%\begin{eqnarray*}
%Z_{f}(\chi,s) = \sum_{\substack{\tau \\
%\textnormal{face of }\Gamma_0}} L_\tau S(\Delta_\tau)(s).
%\end{eqnarray*}
In the last summation also the face $\tau=\Gamma_0$ is included and $S(\Delta_{\Gamma_0})=1$.

If $\Delta_\tau$ is simplicial, say (strictly positively) spanned by primitive linearly independent vectors $\underline{a}_1, \dots, \underline{a}_r \in \Z_{\geq 0}^n$, then
\[ S(\Delta_\tau)(s) = \frac{ \sum_{\underline{h}} q^{\nu(\underline{h}) + N(\underline{h})s} }{\prod_i (q^{\nu(\underline{a}_i) + N(\underline{a}_i)s} - 1)} \]
where the sum runs over $\Z^n \cap \{ \lambda_1 \underline{a}_1 + \dots + \lambda_r \underline{a}_r \, | \, 0 \leq \lambda_i < 1 \}$. In particular if $\mult(\Delta_\tau) = 1$ then the numerator is $1$.
In the non-simplicial case $S(\Delta_\tau)(s)$ is a sum of such expressions (obtained by subdividing $\Delta_\tau$ into simplicial cones).

We clearly see that the real parts of a set of candidate poles
(containing all poles) of the local and global Igusa zeta function are
given by the rational numbers $-\nu(\underline{a})/N(\underline{a})$ for
$\underline{a}$ orthogonal to a facet of the Newton polyhedron at the origin.
Moreover we can restrict to the facets for which the order of $\bar{\chi}$ divides $N(\underline{a})$:
this follows from Lemma~\ref{vanishingofLtau} below.
A fortiori
we can restrict to those for which the order of $\chi$ divides $N(\underline{a})$. We say that such a facet \emph{contributes} a candidate pole to $Z_{f,0}(\chi,s)$ (resp.\ $Z_{f}(\chi,s)$).

We finally remark that if $f$ is nondegenerate over $\C$ w.r.t.\ the compact faces of $\Gamma_0$, then the
couples $(\nu(\underline{a}), N(\underline{a}))$ are part of the numerical data $(\nu_j, N_j)$ associated
to a very explicit (namely, toric) embedded resolution of $f^{-1}\{0\}$ over $F$, that was first described by Varchenko in \cite{V76}.
Thus the fact that we can restrict to the case where the order of $\chi$ divides $N(\underline{a})$ also follows from Igusa's seminal work.

\subsection{The formula of Varchenko for the zeta function of monodromy of $f$ in the origin}

Let $f:(\mathbb{C}^n,0)\rightarrow (\mathbb{C},0)$ be a germ of a holomorphic function. Let $\mathcal{F}$ be the Milnor fibre of the Milnor fibration at the origin associated with $f$ and write $h^i_*:H^i(\mathcal{F},\mathbb{C})\rightarrow H^i(\mathcal{F},\mathbb{C})$, $i\geq 0$, for the monodromy transformations.

The \emph{zeta function of monodromy} at the origin associated to $f$ is \[\zeta_{f,0}(t):=\prod\limits_{i\geq 0}(\det(\id^{i}-th^i_*;H^i(\mathcal{F},\mathbb{C})))^{(-1)^{(i+1)}},\]
where $id^i$ is the identical transformation on $H^i(\mathcal{F},\mathbb{C})$.
One calls $\alpha$ an \emph{eigenvalue of monodromy} of $f$ at the origin if $\alpha$ is an eigenvalue for some $h^i_*:H^i(\mathcal{F},\mathbb{C})\rightarrow H^i(\mathcal{F},\mathbb{C})$.
Denef proved that every eigenvalue of monodromy of $f$ is a zero or a pole of the zeta function of monodromy at some point of $\{f=0\}$ (\cite{D93}).
Varchenko gave in \cite{V76} a formula for $\zeta_{f,0}$ in terms of $\Gamma_0$ if $f$ is nondegenerate w.r.t.\ the compact faces of its Newton polyhedron at the origin $\Gamma_0$.
He defines a function
$\zeta_\tau(t)$ for every compact face $\tau$ of $\Gamma_0$
for which there exists a
subset $I \subset \{1, \ldots , n\} $ with $ \# I =
\textnormal{dim}(\tau)+1$ such that
$\tau \subset L_I:=\{x \in \R^n \, | \, \forall i \not \in I :
x_i = 0\}.$ We will call such faces \emph{V-faces} and we will denote the corresponding index set (resp.\ linear space) to a V-face $\tau$ by $I_\tau$ (resp.\ $L_{I_\tau}$). If a V-face is a simplex, then we will call it a \emph{V-simplex}.

For a face $\tau$ of dimension 0, we put $\textnormal{Vol}(\tau) =
1.$ For every other compact face $\tau$, $\textnormal{Vol}(\tau)$ is defined as the \emph{volume} of $\tau$ for the
volume form $\omega_\tau.$ This is a volume form on
Aff$(\tau),$ the affine space spanned by $\tau,$ such that the
parallelepiped spanned by a lattice-basis of $\Z^n \cap
\textnormal{Aff}(\tau)$ has volume 1. The product
$(\textnormal{dim }\tau)! \textnormal{Vol}(\tau)$ is also called
the \emph{normalized volume} of the face $\tau$ and will be denoted by $\NV(\tau)$. %If $\tau$ is a
%simplicial facet, this normalized volume is equal to the
%multiplicity of the cone spanned by the vertices divided by the
%lattice distance of the facet.

For a V-face $\tau$, let $\sum_{i\in I_\tau}a_i x_i = N(\tau)$ be the equation of
\mbox{Aff}($\tau$) in $L_{I_\tau},$ where $N(\tau)$ and all $a_i$ (for
$i\in I_\tau$) are positive integers and their greatest common divisor is equal to 1. We put
$$\zeta_\tau(t) := \left(1-t^{N(\tau)}\right)^{\NV(\tau)}.$$
In \cite{V76} Varchenko showed that the zeta function
of monodromy of $f$ in the origin is equal to
\begin{eqnarray*}
\zeta_{f,0}(t) = \prod \zeta_{\tau}(t)^{(-1)^{\textnormal{{dim}$(\tau)$}}},
\end{eqnarray*}
where the product runs over all V-faces $\tau$ of $\Gamma_0$.

For a fixed facet $\tau$ of $\Gamma_0$, we say that a V-face $\sigma$ in $\Gamma_0$ contributes w.r.t.\ $\tau$ if $e^{-2 \pi i \nu(\tau)/N(\tau)}$ is a zero of $\zeta_\sigma(t)$.

If $n=3$, the formula of Varchenko for the zeta function of monodromy in the origin has a specific form that we describe below. We first partition every compact facet in simplices, let us
say without introducing new vertices. For such a simplex $\tau$, we define the factor $F_{\tau}$ as in \cite{LVP11}:
\begin{align}
F_{\tau}:=\zeta_{\tau}\prod\limits_{\sigma}\zeta_{\sigma}^{-1}\prod\limits_{p}\zeta_{p},\label{formulaftau}
\end{align}
 where the first product runs over the $1$-dimensional V-faces $\sigma$ in $\tau$ and the second product runs over the $0$-dimensional V-faces $p$ of $\tau$ that are intersection points of two $1$-dimensional V-faces in $\tau$. In \cite[Proposition 8]{LVP11} it is shown that $F_\tau$ is a polynomial.
 Following the formula of Varchenko, the zeta function of monodromy in the origin can be written as
\begin{align}
\zeta_{f,0}(t)=\prod\limits_{\tau}F_{\tau}\prod\limits_{\sigma}\zeta_{\sigma}^{-1}\prod\limits_{p}\zeta_{p},\label{formulan}	
\end{align}
where the first product runs over all $2$-dimensional simplices $\tau$ obtained after subdividing the compact facets and the other products run over $1$-dimensional V-faces $\sigma$ and $0$-dimensional V-faces $p$ for which $\zeta_{\sigma}$, respectively $\zeta_{p}$, was not used in any $F_{\tau}$.

\section{Preliminary results on the normalized volume}

When searching for eigenvalues of monodromy using the formula of Varchenko, one has to compare normalized volumes of compact faces in a facet. This is the motivation for this section.
%
%\begin{pr}\label{nv} Let $\tau$ be a $V$-face and let $\sigma$, $\sigma_1$ and $\sigma_2$ be $V$-faces in $\tau$. Then:
%\begin{enumerate}
%	\item $NV(\tau)=NV(\sigma)\cdot M$, with $M\in \mathbb{N}$.
%	\item $NV(\tau)\cdot NV(\sigma_1\cap\sigma_2)=NV(\sigma_1)\cdot NV(\sigma_2)\cdot M$, with $M\in \mathbb{N}$. If $\sigma_1+\sigma_2=\tau$, then $M=1$ if and only if $N(\sigma_1\cap\sigma_2)=\gcd(N(\sigma_1),N(\sigma_2))$.
%\end{enumerate}
%\end{pr}
For two faces $\sigma$ and $\sigma'$ in a simplicial facet $\tau$, we will denote the smallest face containing $\sigma$ and $\sigma'$ by $\sigma + \sigma'$.
\begin{lm} \label{propcapVfaces}
Let $\sigma$ and $\sigma'$ be two non-disjoint V-faces in a simplicial facet $\tau$. Then $\sigma \cap \sigma'$ and $\sigma + \sigma'$ are also V-faces.
\end{lm}

\noindent \emph{Proof.}
Let $\sigma$ be a $d_1$-dimensional V-simplex and $\sigma'$ a $d_2$-dimensional V-simplex, having $k$ vertices in common.
Suppose that the vertices of $\sigma+\sigma'$ have exactly $s$ zero entries in common. Then one has
\begin{eqnarray*}
s  \leq  n-\#(\sigma+\sigma')   =  n-(d_1+1+d_2+1-k)
\end{eqnarray*}
where (abusing notation) $\#(\sigma + \sigma')$ denotes the number of vertices of $\sigma + \sigma'$.
On the other hand, the vertices of $\sigma \cap \sigma'$ have at most $n-k$ zero entries in common, and so
\begin{eqnarray*}
n-k  \geq  (n-d_1-1)+(n-d_2-1)-s
\end{eqnarray*}
Combining these two inequalities, one finds that they are actually equalities and so $\sigma \cap \sigma'$ and $\sigma + \sigma'$ are V-simplices.
\hfill $\blacksquare$
${}$\\ \\
Recall that for a V-simplex $\tau$, the normalized volume $\NV(\tau)$ is equal to its multiplicity mult$(\tau)$ divided by
its lattice distance $N(\tau)$.
Let $B^{j}(B^j_{1},\ldots,B^j_{n}), 1 \leq j \leq n,$ be the vertices of $\tau$ and let $\sigma$ be a V-face in $\tau$ with vertices $B^1,\ldots,B^k$ and $I_\sigma=\{1,\ldots,k\}$. Then mult$(\tau)$ is the absolute value of the determinant of the matrix
$$\left(\begin{array}{cccccc} B^1_1 & \ldots & B^1_k & 0 & \ldots & 0 \\
\vdots & \ldots & \vdots & \vdots & \ldots & \vdots \\
B^k_1 & \ldots & B^k_k & 0 & \ldots & 0 \\
* & \ldots & * & B^{k+1}_{k+1} & \ldots & B^{k+1}_n \\
\vdots & \ldots & \vdots & \vdots & \ldots & \vdots \\
* & \ldots & * & B^n_{k+1} & \ldots & B^n_n
\end{array}\right).$$
We will denote the matrix
$$M_{\tau,\sigma}:=\left(\begin{array}{ccc}  B^{k+1}_{k+1} & \ldots & B^{k+1}_n \\
 \vdots & \ldots & \vdots \\
B^n_{k+1} & \ldots & B^n_n
\end{array}\right).$$
Then we have that mult$(\tau)=$mult$(\sigma)\left|\det(M_{\tau,\sigma})\right|.$
\begin{pr} \label{propNVdivide}
Let $\tau$ be a simplicial facet of a Newton polyhedron in $\R^n$. If $\sigma$ is a V-face in $\tau$, then $\NV(\sigma) | \NV(\tau)$.
\end{pr}

\noindent \emph{Proof.}
Let us denote the equation of the affine space through $\tau$ resp.\ through $\sigma$ by
$$\mbox{Aff}(\tau) \leftrightarrow a_1x_1 + \ldots + a_nx_n=N(\tau), \qquad \mbox{Aff}(\sigma) \leftrightarrow \frac{a_1x_1 + \ldots + a_kx_k}{\gcd(a_1,\ldots,a_k)}=N(\sigma),$$
with $\gcd(a_1,\ldots,a_n)=1$ and $N(\sigma)=N(\tau)/\gcd(a_1,\ldots,a_k)$.
Let $B^{j}(B^j_{1},\ldots,B^j_{n}), 1 \leq j \leq n$, be the vertices of $\tau$ and $B^{k+1},\ldots,B^n$ the vertices of $\tau$ that are not contained in $\sigma$.
Then we find that
$$ \NV(\tau)=\frac{\NV(\sigma)\left|\mbox{det}(M_{\tau,\sigma})\right|}{\gcd(a_1,\ldots,a_k)}.$$
Let $v_j(B^{k+1}_j,\ldots,B^n_j)^T, k+1 \leq j \leq n,$ be the $j$th column of the matrix $M_{\tau,\sigma}$ and let $\tilde{M_{\tau,\sigma}}$ be the matrix obtained from $M_{\tau,\sigma}$ by replacing the first column by $a_{k+1}v_{k+1}+\ldots + a_nv_n$. For every vertex $B^j$ of $\tau$ we have that
$\gcd(a_1,\ldots,a_k) | a_{k+1}B^j_{k+1} + \ldots + a_nB^j_n$ and hence we find that $\gcd(a_1,\ldots,a_k) | \det(\tilde{M_{\tau,\sigma}})=a_{k+1}\det(M_{\tau,\sigma})$.
Analogously, we obtain that $\gcd(a_1,\ldots,a_k) | a_{j}\det(M_{\tau,\sigma})$, for $k+1 \leq j \leq n$. As we supposed that $\gcd(a_1,\ldots,a_n)=1$, we get that $\gcd(a_1,\ldots,a_k) | \mbox{det}(M_{\tau,\sigma})$, which implies that $\NV(\sigma) | \NV(\tau)$.
\hfill $\blacksquare$

\begin{pr} \label{prNV}
Let $\tau$ be a simplicial facet of a Newton polyhedron in $\R^n$. If $\sigma$ and $\sigma'$ are V-faces in $\tau$ such that $\sigma \cap \sigma' \neq \emptyset$, then
\begin{eqnarray} \label{eqM}
\NV(\tau) \NV(\sigma \cap \sigma') = \NV(\sigma) \NV(\sigma')M, \quad \mbox{for some } M \in \N.
\end{eqnarray}
Moreover, if $\sigma+\sigma'=\tau$, then $M=1$ if and only if $N(\sigma \cap \sigma') = \gcd(N(\sigma),N(\sigma'))$.
\end{pr}

\noindent \emph{Proof.}
As $\sigma+\sigma'$ is also a V-face (see Lemma \ref{propcapVfaces}), it follows by Proposition \ref{propNVdivide} that it is sufficient to prove that
\begin{eqnarray*}
\NV(\sigma+\sigma') \NV(\sigma \cap \sigma') = \NV(\sigma) \NV(\sigma')M, \quad \mbox{for some } M \in \N.
\end{eqnarray*}
Let $B^1,\ldots,B^k,B^{k+1},\ldots,B^r$ be the vertices of $\sigma$ and $B^1,\ldots,B^k,B^{r+1},\ldots,B^s$ be the vertices of $\sigma'$.
Then mult$(\sigma+\sigma')$ is the absolute value of the determinant of the matrix
$$\left(\begin{array}{ccccccccc}
B^1_1  & \ldots & B^1_k  & 0             & \ldots & 0         &  0      & \ldots & 0\\
\vdots & \ldots & \vdots & \vdots        & \ldots & \vdots    & \vdots & \ldots & \vdots \\
B^k_1  & \ldots & B^k_k  & 0             & \ldots & 0         &        0      & \ldots & 0 \\
*      & \ldots & *      & B^{k+1}_{k+1} & \ldots & B^{k+1}_r &     0      & \ldots & 0 \\
\vdots & \ldots & \vdots & \vdots        & \ldots & \vdots    &  \vdots & \ldots & \vdots \\
*      & \ldots & *      & B^r_{k+1}     & \ldots & B^r_{r}   &  0      & \ldots & 0 \\
*      & \ldots & *      & 0             & \ldots & 0         & B^{r+1}_{r+1} & \ldots & B^{r+1}_s\\
\vdots & \ldots & \vdots & \vdots        & \ldots & \vdots    & \vdots & \ldots & \vdots  \\
*      & \ldots & *      & 0             & \ldots & 0         &  B^s_{r+1}   & \ldots & B^{s}_s
\end{array}\right).$$
%We call
%%M_{\sigma \cap \sigma'}:=\left(\begin{array}{ccc}
%%B^1_1  & \ldots & B^1_k  \\
%%\vdots & \ldots & \vdots \\
%%B^k_1  & \ldots & B^k_k  \end{array}\right), \quad
%%M_{\sigma,\sigma \cap \sigma'}:=\left(\begin{array}{ccc} B^{k+1}_{k+1} & \ldots & B^{k+1}_r  \\
%% \vdots        & \ldots & \vdots   \\
%% B^r_{k+1}     & \ldots & B^r_{r} \end{array}\right),$$
%% $$M_{\sigma',\sigma \cap \sigma'}:=\left(\begin{array}{ccc}
%% B^{r+1}_{n-s+r+1} & \ldots & B^{r+1}_n\\
%% \vdots & \ldots & \vdots  \\
%% B^s_{n-s+r+1}   & \ldots & B^{s}_n
%% \end{array}\right), \quad
%$$M_{\tau,\sigma+\sigma'}:=\left(\begin{array}{ccc}
% B^{s+1}_{r+1} & \ldots & B^{s+1}_{n-s+r} \\
% \vdots & \ldots & \vdots  \\
% B^n_{r+1}     & \ldots & B^{n}_{n-s+r}
% \end{array}\right).$$
We write
\begin{eqnarray*}
\mbox{Aff}(\sigma+\sigma')  \leftrightarrow  a_1x_1 + \ldots + a_sx_s=N(\sigma+\sigma'),  \quad \mbox{with } \gcd(a_1,\ldots,a_s)=1, \\
\alpha:=\gcd(a_1,\ldots,a_k), \beta:=\gcd(a_{k+1},\ldots,a_r) \mbox{ and } \gamma:=\gcd(a_{r+1},\ldots,a_s).
\end{eqnarray*}
Then we have
\begin{eqnarray*}
\mbox{Aff}(\sigma) & \leftrightarrow & \frac{a_1x_1 + \ldots + a_kx_k+a_{k+1}x_{k+1}+ \ldots + a_rx_r}{\gcd(\alpha,\beta)}=\frac{N(\sigma+\sigma')}{\gcd(\alpha,\beta)}=N(\sigma), \\
\mbox{Aff}(\sigma') & \leftrightarrow & \frac{a_1x_1 + \ldots + a_kx_k+a_{r+1}x_{r+1}+ \ldots + a_sx_s}{\gcd(\alpha, \gamma)}=\frac{N(\sigma+\sigma')}{\gcd(\alpha,\gamma)}=N(\sigma'),\\
\mbox{Aff}(\sigma \cap \sigma') & \leftrightarrow & \frac{a_1x_1 + \ldots + a_kx_k}{\alpha}=\frac{N(\sigma+\sigma')}{\alpha}=N(\sigma \cap \sigma').
\end{eqnarray*}
%If we denote $\delta:=\det(M_{\tau,\sigma+\sigma'})$, then
By using Proposition \ref{propNVdivide} we get
\begin{eqnarray*}
\NV(\sigma+\sigma') \NV(\sigma \cap \sigma') = \NV(\sigma) \NV(\sigma')\frac{\alpha}{\gcd(\alpha,\beta)\gcd(\alpha,\gamma)}.
\end{eqnarray*}
As gcd$(\alpha,\beta,\gamma)=1$, the quotient $\alpha/(\gcd(\alpha,\beta)\gcd(\alpha,\gamma))$ is an integer.
%We will now prove that
%\begin{eqnarray} \label{eqp}
%M:=\frac{\delta\alpha}{\gcd(\alpha,\beta)\gcd(\alpha,\gamma)}=\frac{\delta lcm(\alpha,\beta)lcm(\alpha,\gamma)}{\alpha\beta\gamma} \in \N.
%\end{eqnarray}
%Let $p$ be a prime number. Let $g$ (resp.\ $h$, resp.\ $i$, resp.\ $j$) be the maximal power of $p$ such that $p^g | \alpha$ (resp.\ $p^h | \beta$, resp.\ $p^i | \gamma$, resp.\ $p^j | \delta$). To prove (\ref{eqp}), we will show that
%\begin{eqnarray}
%g+h+i \leq j + \max\{g,h\}+\max\{g,i\}.
%\end{eqnarray}
%If $g,h$ or $i$ is equal to $0$, then this inequality follows directly.
%If $g,h$ and $i$ are different from $0$, then let $l:=\min\{g,h,i\}$.
%As $$g+h+i \leq \min\{g,h,i\} + \max\{g,h\} + \max\{g,i\},$$
%it is now sufficient to prove that $l \leq j$ or $p^l | \delta$.
%The equation of $\mbox{Aff}(\tau)$ modulo $p^l$ becomes
%\begin{eqnarray*}
%a_{s+1}x_{s+1} + \ldots a_n x_n \equiv 0 \mbox{ mod } p^l.
%\end{eqnarray*}
%As $\gcd(a_1,\ldots,a_n)=1$, there is an $a_t$, $s+1 \leq t \leq n$ such that $p$ does not divide $a_t$ and hence
%$a_t$ is a unit in $\Z_{/{(p^l)}}$. This implies that the vector $(B^{s+1}_t,\ldots,B^n_t)^T$ is a linear combination of the other column vectors in $M_{\tau,\sigma+\sigma'}$ in $\Z_{/{(p^l)}}$ and hence that $p^l | \delta$.

To prove the second statement, let $\sigma$ and $\sigma'$ be two V-faces in a simplicial facet $\tau$ such that $\sigma+\sigma'=\tau$.
Then one easily shows that $N(\tau)=\mbox{lcm}(N(\sigma),N(\sigma'))$ and one can then write
$$M=\frac{\alpha}{\gcd(\alpha,\beta)\gcd(\alpha,\gamma)}=\frac{\gcd(N(\sigma),N(\sigma'))}{N(\sigma \cap \sigma')}.$$
\hfill $\blacksquare$
\begin{co}
Let $\sigma$ and $\sigma'$ be two V-faces in a simplicial facet $\tau$. If $\sigma$ and $\sigma'$ contribute w.r.t.\ $\tau$ and if $\sigma \cap \sigma'$ does not, then $M \geq 2$ in Equation \emph{(\ref{eqM})}.
\end{co}

\section{Some character sums} \label{sectioncharsums}

In order to prove the holomorphy conjecture, we will have to show that some candidate poles of $Z_{f,0}(\chi,s)$ (resp.\ $Z_{f}(\chi,s)$) are false poles.
These proofs rely on the computation of some character sums. We first recall some well-known properties of character sums over finite fields which we will need when treating $B_1$-facets. We then study a specific character sum (see Proposition \ref{eqchar2}) which shows up when proving fakeness of some other candidate pole.

\begin{lm} \label{vanishingofLtau}
Let $a_1, \dots, a_n, N \in \Z$ and let $\chi$ be a multiplicative character
of $\F_q^\times$ whose order is not a divisor of $N$. Let $f \in \F_q[x_1, \dots, x_n]$ be such that
each exponent $(k_1, \dots, k_n)$ appearing in $f$ satisfies $a_1k_1 + \dots + a_n k_n = N$. Then
\[ \sum_{(x_1,\ldots,x_n) \in (\F_q^\times)^n} \chi(f(x_1, \dots, x_n)) = 0. \]
\end{lm}
\noindent \emph{Proof.}
Pick $u \in \F_q^\times$ such that $\chi(u^N) \neq 1$. Then the left hand side equals
$$ \sum_{(x_1,\ldots,x_n) \in (\F_q^\times)^n}
\chi(f(u^{a_1}x_1, \dots, u^{a_n}x_n)) = \chi\left(u^N\right)\sum_{(x_1,\ldots,x_n) \in (\F_q^\times)^n} \chi(f(x_1, \dots, x_n))$$
from which the property follows.
\hfill $\blacksquare$

\begin{lm} \label{eqchar0}
Let $a \in \N$ and let $\chi$ be a multiplicative character of $\F_q^\times$ whose order is not a divisor of $a$, then
$\sum_{x \in \F_q^\times} \chi(x^a)=0.$
\end{lm}
\noindent \emph{Proof.}
Take $f(x) = x^a$ in the foregoing lemma.
\hfill $\blacksquare$

\begin{lm} \label{eqchar1}
Let $f$ be a polynomial and $g$ be a monomial (possibly equipped with a non-zero coefficient) over $\F_q$ in the variables $x_2, \ldots, x_n$, and let $\chi$ be a non-trivial multiplicative character of $\F_q^\times$. Then
$$\sum_{(x_1,\ldots,x_n) \in (\F_q^\times)^n} \chi(f(x_2,\ldots,x_n)+x_1g(x_2,\ldots,x_n))=-\sum_{(x_2,\ldots,x_n) \in (\F_q^\times)^{n-1}}\chi(f(x_2,\ldots,x_n)).$$
\end{lm}
\noindent \emph{Proof.}
One can write
\begin{eqnarray*}
& & \sum_{(x_1,\ldots,x_n) \in (\F_q^\times)^n} \chi(f(x_2,\ldots,x_n)+x_1g(x_2,\ldots,x_n)) \\ & = & \sum_{(x_2,\ldots,x_n) \in (\F_q^\times)^{n-1}} \sum_{x_1 \in \F_q^\times}\chi(f(x_2,\ldots,x_n)+x_1g(x_2,\ldots,x_n)) \\
 & =& \sum_{(x_2,\ldots,x_n) \in (\F_q^\times)^{n-1}}\left(\sum_{u \in \F_q}\chi(u) - \chi(f(x_2,\ldots,x_n))\right) \\
 & = & -\sum_{(x_2,\ldots,x_n) \in (\F_q^\times)^{n-1}}\chi(f(x_2,\ldots,x_n)),
\end{eqnarray*}
where we used Lemma \ref{eqchar0} in the last step.
\hfill $\blacksquare$

\begin{pr} \label{eqchar2}
Let $\chi$ be a multiplicative character of $\F_q^\times$ such that its order does not divide $a \in \N$. Let $\alpha \in \F_q$ and $\beta, \gamma \in \F_q^\times$.
Then
$$\sum_{ x,y,z \in \mathbb{F}_q^\times } \chi(\alpha x^a + \beta x^i y^2 + \gamma x^i z^2) \, = \, - \sum_{ x,y \in \mathbb{F}_q^\times } \chi( \alpha x^a + \beta x^i y^2) \,
- \sum_{ x,z \in \mathbb{F}_q^\times } \chi( \alpha x^a + \gamma x^i z^2), $$
with $i \in \N$.
\end{pr}
\noindent \emph{Proof.}
First notice that we can reduce to the cases $i=0$ or $i=1$. We can also assume that $q$ is odd: indeed if $q$ is even then we can
replace $y^2$ by $y$ and $z^2$ by $z$, from which one sees that all sums are zero, for instance by using Lemma~\ref{eqchar0} and Lemma~\ref{eqchar1}.
Let
\[ \varepsilon = \left\{ \begin{array}{ll} 2 & \text{if $-\beta/\gamma$ is a square in $\mathbb{F}_q$,} \\ 0 & \text{if not.} \\ \end{array} \right. \]
If $i=0$, then for each $c \in \mathbb{F}_q^\times$ define
\[ L_c :=  \# \{ \, (x,y,z) \in \mathbb{F}_q^{\times 3} \,| \, \alpha x^a + \beta y^2 + \gamma z^2 = c \, \}, \qquad M_c := \# \{ \, x \in \mathbb{F}_q^\times \,| \, \alpha x^a  = c \, \}, \]
\[ N_{1,c} :=  \# \{ \, (x,y) \in \mathbb{F}_q^{\times 2} \,| \, \alpha x^a + \beta y^2 = c \, \},
\qquad N_{2,c} :=  \# \{ \, (x,z) \in \mathbb{F}_q^{\times 2} \,| \, \alpha x^a + \gamma z^2 = c \, \}. \]
We rewrite the first equation as
\begin{equation} \label{mainequation}
  \beta y^2 + \gamma z^2 = c - \alpha x^a.
\end{equation}
For each value of $x \in \mathbb{F}_q^\times$ this defines a conic in the variables $y$ and $z$.
In the $M_c$ cases where $c - \alpha x^a = 0$ the conic carries $\varepsilon (q-1) + 1$ points $(y,z) \in \mathbb{F}_q^2$.
If $c - \alpha x^a \neq 0$ then Equation (\ref{mainequation}) defines a smooth conic having $\varepsilon$ points at infinity. As every projective nonsingular curve of genus $0$ over a finite field $\F_q$ has always $q+1$ points (see \cite{W}), the conic carries $q+1 - \varepsilon$ points in $\mathbb{F}_q^2$.
We conclude that there are
\[ (\varepsilon (q-1) + 1) M_c + (q + 1 - \varepsilon) (q-1 - M_c)\]
solutions $(x,y,z) \in \mathbb{F}_q^\times \times \mathbb{F}_q^2$ to Equation (\ref{mainequation}). Because in $L_c$ there
are $M_c$ points of the form $(x,0,0)$, $N_{1,c}$ points of the form $(x,y,0)$ with $y \neq 0$, and
$N_{2,c}$ points of the form $(x,0,z)$ with $z \neq 0$,
%\begin{itemize}
%  \item $N_c$ points of the form $(x,0,0)$,
%  \item $M_c$ points of the form $(x,y,0)$ with $y \neq 0$,
%  \item $M_c$ points of the form $(x,0,z)$ with $z \neq 0$,
%\end{itemize}
we conclude that
\[  L_c  =  (\varepsilon (q-1) + 1) M_c + (q + 1 - \varepsilon) (q-1 - M_c) - M_c - N_{1,c} - N_{2,c}. \]
Then for some constant $\lambda$, it holds that $L_c = - N_{1,c} - N_{2,c} + \lambda M_c$.
Now note that
\[ S_1 := \sum_{ x,y,z \in \mathbb{F}_q^\times } \chi(\alpha x^a + \beta y^2 + \gamma z^2) = \sum_{c \in \mathbb{F}_q^\times} L_c \chi(c),\]  \[ S_2 := \sum_{ x,y \in \mathbb{F}_q^\times } \chi(\alpha x^a + \beta y^2) = \sum_{c \in \mathbb{F}_q^\times} N_{1,c} \chi(c), \ \  S_3 := \sum_{ x,z \in \mathbb{F}_q^\times } \chi(\alpha x^a + \gamma z^2) = \sum_{c \in \mathbb{F}_q^\times} N_{2,c} \chi(c), \]
\[ 0 = \chi(\alpha) \sum_{ x \in \mathbb{F}_q^\times } \chi(x^a) = \sum_{ x \in \mathbb{F}_q^\times } \chi(\alpha x^a) = \sum_{c \in \mathbb{F}_q^\times} M_c \chi(c). \]
In the last case, the first equality follows by Lemma \ref{eqchar0}.
Plugging in the expression for $L_c$ in $S_1$ we find
\[ S_1 = - \sum_{c \in \mathbb{F}_q^\times}N_{1,c} \chi(c)  - \sum_{c \in \mathbb{F}_q^\times}N_{2,c} \chi(c) + \lambda \sum_{c \in \mathbb{F}_q^\times}N_c \chi(c) = - S_2 - S_3  \]

If $i=1$, one instead of (\ref{mainequation}) uses the conic
\begin{equation*}
  \beta y^2 + \gamma z^2 = \frac{c - \alpha x^a}{x}
\end{equation*}
and proceeds analogously.
\hfill $\blacksquare$\\

Note that the exponents $(a,0,0)$, $(i,2,0)$, $(i,0,2)$ are contained in the hyperplane $2k_1 + (a-i)k_2 + (a-i)k_3 = 2a$, so under the stronger assumption
that the order of $\chi$ does not divide $2a$, or under the additional assumption that $a-i$ is even, we see from Lemma~\ref{vanishingofLtau} that
all sums in the statement of the proposition are actually zero.

\section{A proof of the holomorphy conjecture for nondegenerate surface singularities}

Let $f(\underline{x})$ be as in Subsection~\ref{subsection_nondegsing} and
assume that it is nondegenerate over $\C$ with respect to the compact faces (resp.\ the faces) of its Newton polyhedron
at the origin $\Gamma_0$. Let $K$ be a non-archimedean completion with valuation ring $R$ and residue field $\F_q$, such
that $\bar{f}$ is nondegenerate over $\F_q$ with respect to the compact faces (resp.\ the faces) of $\Gamma_0$. Let $\chi : R^\times \rightarrow \mathbb{C}^\times$ be a non-trivial character of conductor $1$.
If $Z_{f,0}(\chi,s)$ (resp.\ $Z_{f}(\chi,s)$) is not holomorphic on $\mathbb{C}$, then by the material from Subsection~\ref{subsigchar} it has a pole with real part equal to $-\nu(\tau) / N(\tau)$ for some facet $\tau$ of $\Gamma_0$ for which the order of $\bar{\chi}$
divides $N(\tau)$. Here as before $\bar{\chi}$ denotes the unique character of $\F_q^\times$ associated to $\chi$.

For some facets $\tau$, in particular the $B_1$-facets and the $X_2$-facets which we introduce here, we will mostly have to prove that $-\nu(\tau)/N(\tau)$ cannot be the real part of a pole of $Z_{f,0}(\chi,s)$ (resp.\ $Z_{f}(\chi,s)$). For the other facets, we will prove that $e^{-2\pi i/N(\tau)}$ is an eigenvalue of monodromy of $f$ at some point of $f^{-1} \{ 0 \}$ and we will thus obtain that the order of $\chi$ (which we recall divides the order of $\bar{\chi}$) divides the order of some eigenvalue of monodromy at some point of $f^{-1}\{ 0 \}$.

Let us first recall the notion of $B_1$-facets, introduced in \cite{LVP11}.
A simplicial facet of an $n$-dimensional Newton polyhedron ($n\geq
2$) is a \emph{$B_1$-simplex w.r.t.\ the variable $x_i$}
if it is a simplex with $n-1$ vertices in the coordinate
hyperplane $x_i=0$ and one vertex at distance one of this
hyperplane.
We say that a facet $\tau$ of an $n$-dimensional Newton polyhedron is \emph{non-compact for
the variable $x_j$} ($1 \leq j \leq n$) if for every point $p \in
\tau$ the point $p + (0, \ldots , 0,1,0, \ldots,0) \in \tau$,
where $(0, \ldots , 0,1,0, \ldots,0)$ is an $n$-tuple with 1 at
place $j$ and 0 everywhere else. We define the maps
$\pi_j: \R^n \to \R^{n-1}: (x_1, \ldots , x_n) \mapsto (x_1,
\ldots, \widehat{x_j}, \ldots, x_n)$ for $j=1, \ldots , n$.
A non-compact facet $\tau$ of an $n$-dimensional Newton polyhedron
($n\geq 3$) is a (non-compact) \emph{$B_1$-facet w.r.t.\
the variable $x_i$} if $\tau$ is non-compact for exactly one
variable $x_j$ and if $\pi_j(\tau)$ is a $B_1$-simplex in
$\R^{n-1}$ w.r.t.\ $x_i$.
A \emph{$B_1$-facet} is a $B_1$-simplex or a non-compact
$B_1$-facet w.r.t.\ some variable.

\begin{de}
A \emph{facet of type $X_2$} in a $3$-dimensional Newton polyhedron is a facet whose vertices (up to permutation of the coordinates) are of the form
$p(a,0,0), q(x_1,0,2),$ $r(x_2,2,0)$ with $a - x_2$ and $a - x_1$ both odd. %and $x_1 \neq 0$.
\end{de}

\subsection{Determination of a set of eigenvalues} \label{subseigenvalues}

In \cite[Theorem 10]{LVP11} Van Proeyen and the third author proved that $e^{-2\pi i\nu(\tau)/N(\tau)}$ is an eigenvalue of monodromy of $f$ at some point of $f^{-1}\left\{0\right\}$, whenever $\tau$ is not a $B_1$-facet. We will now show that $e^{-2\pi i/N(\tau)}$ is also an eigenvalue of monodromy at some point of $f^{-1}\left\{0\right\}$, except for some cases which will be treated in Subsection \ref{subsB1}.
We here rely on Proposition \ref{prNV}
to get a more conceptual proof.

We first divide the compact facets in simplices $\tau$, without introducing new vertices.

\begin{pr}\label{prdeler}
%Let $\tau$ be a compact simplicial facet of a $3$-dimensional Newton polyhedron $\Gamma_0$ which is
%. If $-\nu(\tau)/N(\tau)$ is the real part of a pole of $Z_{f,0}(\chi,s)$ (resp.\ $Z_{f}(\chi,s)$) and if $\tau$ is
If $\tau$ is not of type $B_1$ nor of type $X_2$, then $e^{-2\pi i /N(\tau)}$ is a zero of $F_{\tau}$.
\end{pr}

\noindent \emph{Proof.}
\textsc{Case 1: $\tau$ does not contain a segment in a coordinate plane}.\\
By formula (\ref{formulaftau}), $F_{\tau}=\zeta_{\tau}=\left(1-t^{N(\tau)}\right)^{\NV(\tau)}$ and $e^{-2\pi i/N(\tau)}$ clearly is a zero of $F_{\tau}$.
\\ \\
\textsc{Case 2: $\tau$ contains exactly one $1$-dimensional V-face $\sigma$}.
\\
In this case, we have \[F_{\tau}=\frac{\zeta_{\tau}}{\zeta_{\sigma}}=\frac{\left(1-t^{N(\tau)}\right)^{\NV(\tau)}}{\left(1-t^{N(\sigma)}\right)^{\NV(\sigma)}}.\]
Then $e^{-2\pi i/N(\tau)}$ is a zero of $F_{\tau}$ unless $N(\sigma)=N(\tau)$ and $\NV(\sigma)= \NV(\tau)$. One easily checks that then $\tau$ would be a $B_1$-facet.
\\ \\
\textsc{Case 3: $\tau$ contains exactly two $1$-dimensional V-faces $\sigma_1$ and $\sigma_2$.}
\\
In this situation, \[F_{\tau}=\frac{\zeta_{\tau}\zeta_p}{\zeta_{\sigma_1} \zeta_{\sigma_2}}=\frac{(1-t^l)\left(1-t^{N(\tau)}\right)^{\NV(\tau)}}{\left(1-t^{N(\sigma_1)}\right)^{\NV(\sigma_1)}\left(1-t^{N(\sigma_2)}\right)^{\NV(\sigma_2)}},\] where w.l.o.g.\ $\left\{p\left(l,0,0\right)\right\}=\sigma_1\cap\sigma_2$.

If $N(\sigma_1) \neq N(\tau)$ or $N(\sigma_2) \neq N(\tau)$, then see Case 1 and Case 2.
If $N(\tau)= N(\sigma_1)= N(\sigma_2)$, then $F_{\tau}=(1-t^l)\left(1-t^{N(\tau)}\right)^{\NV(\tau)-\NV(\sigma_1)-\NV(\sigma_2)}$.

\medskip
\textit{Case 3.1}: If $N(p)=N(\tau)$, then by Proposition \ref{prNV}, $\NV(\tau)=\NV(\sigma_1)\NV(\sigma_2)$ and hence $F_{\tau}=\left(1-t^{N(\tau)}\right)^{\left(\NV(\sigma_1)-1\right)\left(\NV(\sigma_2)-1\right)}$.
If $\NV(\sigma_1)$ or $\NV(\sigma_2)$ would be equal to $1$, then it would result that $\NV(\tau)=\NV(\sigma_i)$, for some $i \in \{1,2\}$ and again $\tau$ would be a $B_1$-facet. Consequently, $e^{-2\pi i/N(\tau)}$ is a zero of $F_{\tau}$.

\medskip
\textit{Case 3.2}: Suppose that $N(p)\neq N(\tau)$.
 By Proposition \ref{prNV}, we have $\NV(\tau)=M \NV(\sigma_1) \NV(\sigma_2)$, with $M \geq 2$. If not both $\NV(\sigma_1)$ and $\NV(\sigma_2)$ are equal to $1$, then one easily deduces that $\NV(\tau)-\NV(\sigma_1)-\NV(\sigma_2)>0$.
If $\NV(\sigma_1)=\NV(\sigma_2)=1$, then $\NV(\tau)-\NV(\sigma_1)-\NV(\sigma_2)>0$ if and only if $M > 2$.
It remains thus to study the case $\NV(\tau)=M=2$, $\NV(\sigma_1)=\NV(\sigma_2)=1$. As we supposed that $N(\tau)= N(\sigma_1)= N(\sigma_2)$, the vertices of $\tau$ are then $p(N(\tau)/2,0,0), q(x_1,0,2), r(x_2,2,0)$, and
$$\mbox{Aff}(\tau) \leftrightarrow 2x + (N(\tau)/2 - x_2)y + (N(\tau)/2 - x_1)z = N(\tau).$$
From $N(\sigma_1)=N(\sigma_2)=N(\tau)$ it follows that $N(\tau)/2 - x_2$ and $N(\tau)/2 - x_1$ are odd and hence $\tau$ is of type $X_2$.
\\ \\
\textsc{Case 4: $\tau$ contains three $1$-dimensional V-faces $\sigma_1, \sigma_2$ and $\sigma_3$.}\\
In this situation  \[F_{\tau}=\frac{\zeta_{\tau} \zeta_{p} \zeta_{q} \zeta_{r}}{\zeta_{\sigma_1} \zeta_{\sigma_2} \zeta_{\sigma_3}},\]
with $p=\sigma_1 \cap \sigma_2$, $q=\sigma_1 \cap \sigma_3$ and $r=\sigma_2 \cap \sigma_3$. We suppose that $N(\tau)=N(\sigma_1)=N(\sigma_2)=N(\sigma_3)$, if not then we fall back on one of the previous cases.

\medskip
\textit{Case 4.1}: If $N(\tau)= N(\sigma_1\cap\sigma_2)=N(\sigma_1\cap\sigma_3)=N(\sigma_2\cap\sigma_3)$, then, by Proposition \ref{prNV}, $\NV(\tau)=\NV(\sigma_1)\NV(\sigma_2)=\NV(\sigma_1)\NV(\sigma_3)=\NV(\sigma_2)\NV(\sigma_3)$, and thus $\NV(\sigma_1)=\NV(\sigma_2)=\NV(\sigma_3)$. Then $F_\tau$ becomes $$F_{\tau}=\frac{\left(1-t^{N(\tau)}\right)^{\NV(\sigma_1)^2+3}}{\left(1-t^{N(\tau)}\right)^{3\NV(\sigma_1)}}.$$ Since $\NV(\sigma_1)^2+3 >3\NV(\sigma_1)$, it follows that $e^{-\frac{2\pi i}{N(\tau)}}$ is a zero of $F_{\tau}$.

\medskip
\textit{Case 4.2}: If $N(\tau)= N(\sigma_1\cap\sigma_2)=N(\sigma_2\cap\sigma_3)\neq N(\sigma_1\cap\sigma_3)$, then Proposition \ref{prNV} yields $\NV(\tau)=\NV(\sigma_1)\NV(\sigma_2)=\NV(\sigma_2)\NV(\sigma_3)=M \NV(\sigma_1)\NV(\sigma_3)$, with $M \geq 2$. We thus get $\NV(\sigma_3)=\NV(\sigma_1)$ and $\NV(\sigma_2)=M \NV(\sigma_1)$ and we find then
$$\mbox{Aff}(\tau) \leftrightarrow x+My+z=N(\tau),$$
with $p(N(\tau),0,0),q(0,N(\tau)/M,0)$ and $r(0,0,N(\tau))$.
In this case $e^{-\frac{2\pi i}{N(\tau)}}$ is a zero of $F_{\tau}$ if and only if $\NV(\tau)+2>\NV(\sigma_1)+\NV(\sigma_2)+\NV(\sigma_3)$, or equivalently, if $\left(M \NV(\sigma_1)-2\right) \left(\NV(\sigma_1)-1\right)>0$. This is always the case, as $\NV(\sigma_1)=N(\tau)/M=1$ would imply that $\tau$ is a $B_1$-facet.

\medskip
\textit{Case 4.3}: If $N(\tau)= N(\sigma_1\cap\sigma_2)$, $N(\tau)\neq N(\sigma_1\cap\sigma_3)$ and $N(\tau)\neq N(\sigma_2\cap\sigma_3)$, then by Proposition \ref{prNV} one has $\NV(\tau)=\NV(\sigma_1)\NV(\sigma_2)=M_1 \NV(\sigma_1)\NV(\sigma_3)=M_2 \NV(\sigma_2)\NV(\sigma_3)$, with $M_1 \geq 2$ and $M_2 \geq 2$.
In this configuration we have
$$\mbox{Aff}(\tau) \leftrightarrow x+ky+lz=N(\tau),$$
$p(N(\tau),0,0),q(0,N(\tau)/k,0)$ and $r(0,0,N(\tau)/l)$ with $\gcd(k,l) = 1$.
Then we find that $M_1=k$, $M_2=l$ and hence $\NV(\sigma_2)=k\NV(\sigma_1)/l$ and $\NV(\sigma_3)=\NV(\sigma_1)/l$.
%But $NV(\tau)=g_{\tau}=am$, $NV(\sigma_1)=a=\\gcd(a,l)$, $NV(\sigma_2)=m=\\gcd(m,l)$, $NV(\sigma_3)=\frac{am}{l}$. Thus, since $a|l$ and $m|l$, we get $l=ak_1=mk_2$, with $k_1,k_2\in \mathbb{N}_{\geq 2}$. Moreover, since $k_1|m$ and $k_2|a$, it follows that $m=uk_1$, $a=vk_2$, with $u,v\in \mathbb{N}\backslash 0$. Hence, because $u$ and $v$ should be equal, $m=k_1u$, $a=k_2u$, $l=k_1k_2u$ and $NV(\tau)=k_1k_2u^2$, $NV(\sigma_1)=k_2u$, $NV(\sigma_2)=k_1u$, $NV(\sigma_3)=u$.
In this case, $e^{-2\pi i/N(\tau)}$ would be a zero of $F_{\tau}$ if and only if $\NV(\tau)+1 > \NV(\sigma_1)+\NV(\sigma_2)+\NV(\sigma_3)$, or equivalently,
$$k\NV(\sigma_1)^2 - (k + l + 1)\NV(\sigma_1) + l > 0.$$
This is true because $\NV(\sigma_1) \geq l$ while the largest real root of the polynomial on the left hand side is
\[ \frac{k+l+1 + \sqrt{(k + l + 1)^2 - 4kl}}{2k} < l;\]
the latter inequality holds because one easily rewrites it as $kl > k + l$, which holds since $k,l \geq 2$ and $k=l=2$ is excluded
by coprimality.

\medskip
\textit{Case 4.4}: If $N(\tau)\neq N(\sigma_1\cap\sigma_2)$, $N(\tau)\neq N(\sigma_1\cap\sigma_3)$ and $N(\tau)\neq N(\sigma_2\cap\sigma_3)$,
 then by Proposition \ref{prNV} one has $\NV(\tau)=M_1\NV(\sigma_1)\NV(\sigma_2)=M_2 \NV(\sigma_1)\NV(\sigma_3)=M_3 \NV(\sigma_2)\NV(\sigma_3)$, with $M_1 \geq 2$, $M_2 \geq 2$ and $M_3 \geq 2$.
In this configuration we have
$$\mbox{Aff}(\tau) \leftrightarrow kx+ly+mz=N(\tau),$$
$p(N(\tau)/k,0,0),q(0,N(\tau)/l,0)$ and $r(0,0,N(\tau)/m)$ with $k,l,m$ pairwise coprime. Then we find
that $M_1 = k$, $M_2 = l$, $M_3 = m$ and hence $\NV(\sigma_2) = l \NV(\sigma_1)/m$ and $\NV(\sigma_3) = k \NV(\sigma_1)/m$.
In this case we want to establish that
$\NV(\tau) > \NV(\sigma_1)+ \NV(\sigma_2)+ \NV(\sigma_3)$, or equivalently that
$klm \NV(\sigma_1) > k + l + m$. This follows from $\NV(\sigma_1) \geq 1$
and $klm \geq 4 \max \{k,l,m\} > 3 \max \{k,l,m\} \geq k + l + m$. \hfill $\blacksquare$

%Let's suppose that $NV(\sigma_1)\geq NV(\sigma_2)\geq NV(\sigma_3)$. Proposition \ref{prNV} implies that $NV(\tau)=M NV(\sigma_1)NV(\sigma_2)$
%with $M \geq 2$. We will prove that
%\begin{eqnarray} \label{eq44}
%M NV(\sigma_1)NV(\sigma_2)>NV(\sigma_1)+2NV(\sigma_2),
%\end{eqnarray}
%which then implies the needed inequality. \\Let $p(N(\tau)/k,0,0),q(0,N(\tau)/l,0)$ and $r(0,0,N(\tau)/m)$.
%\medskip
%\textit{Case 4.4.1}: We suppose that $M\geq 3$. If $NV(\sigma_2)\geq 2$, then $M NV(\sigma_2)-1 > 2 NV(\sigma_2)$ and this implies (\ref{eq44}).
%If $NV(\sigma_2)=1$, then also $NV(\sigma_3)=1$. If $NV(\sigma_1)\geq 2$, then (\ref{eq44}) follows immediately. If also $NV(\sigma_1)=1$, then
%$N(\tau)=kl=km=lm$ and hence $k=l=m$. However, if $N(\sigma_i)=N(\tau)$ ($1 \leq i \leq 3$), then $k=l=m=1$ but then $\tau$ would fall under
%Case 4.1.

%\medskip
%\textit{Case 4.4.2}: We suppose that $M=2$. If $NV(\sigma_1)\geq 3$, then $NV(\sigma_1) < 2(NV(\sigma_1)-1)$, which implies (\ref{eq44}). If %
%$NV(\sigma_1)=2$ and $NV(\sigma_2)\geq 2$, then again the inequality (\ref{eq44}) clearly is satisfied. If $NV(\sigma_1)=2$ and
%$NV(\sigma_2)=1$, then we have $NV(\sigma_3)=1$ and $N(\tau)=2kl=km=lm$ and hence $k=l=m/2$.
%Again, we find that $\tau$ does not fall under Case 4.4.
%If $NV(\sigma_1)=1$, then also $NV(\sigma_2)=NV(\sigma_3)=1$. As $M=2$, it follows that $NV(\tau) \neq 3$, hence $NV(\tau) \neq
%NV(\sigma_1)+NV(\sigma_2)+NV(\sigma_3)$ and we can conclude that $e^{-2\pi i/N(\tau)}$ is a zero of $F_{\tau}$.

\subsection{On false poles contributed by $B_1$-facets and $X_2$-facets}\label{subsB1}

In \cite[Proposition 9.6]{BV15} it is shown that if a candidate pole contributed only by $B_1$-facets is an actual pole of $Z_{f,0}(\chi,s)$, then it is contributed by two $B_1$-facets w.r.t.\ different variables having a $1$-dimensional intersection.
We will here show that even in that situation the candidate pole is almost always a false pole of $Z_{f,0}(\chi,s)$. We need this precision here, because for the holomorphy conjecture one has to verify if $1/N_j$ gives rise to an eigenvalue of monodromy, rather than the quotient $\nu_j/N_j$ (that might be simplifiable). We also study when candidate poles of $Z_{f}(\chi,s)$ corresponding to $B_1$-facets are false poles.
Finally we provide some facets of type $X_2$ that give rise to fake poles of $Z_{f,0}(\chi,s)$ and $Z_{f}(\chi,s)$.

We again assume that the compact facets have been subdivided into simplices without introducing new vertices; this guarantees that every vertex is equipped with a non-zero coefficient.
Then towards our study of contributions of configurations of $B_1$-facets, we make the following observations (which hold up to permutation of the coordinates).
\\ \\
\textsc{Fact 1}: A vertex $P(1,\cdot,\cdot)$ does not contribute. Indeed,
as $\chi$ is not the trivial character (and so $\bar{\chi}$ neither is trivial), one immediately deduces from Lemma \ref{eqchar0} that the contribution of $P$ is equal to $0$.
\\ \\
\textsc{Fact 2}: A vertex $P(a,0,0)$ does not contribute if the order of $\bar{\chi}$ is not a divisor of $a$ (by Lemma \ref{eqchar0}).
\\ \\
\textsc{Fact 3}: A segment $\sigma:=PQ$ with $P(1,1,b)$ and $Q(0,0,a)$ does not contribute if the order of $\bar{\chi}$ is not a divisor of $a$. To compute the contribution of $\sigma$, we consider
\begin{eqnarray*}
L_\sigma=q^{-3}\sum_{(x,y,z) \in (\F_q^\times)^3} \bar{\chi}(c_{0,0,a}z^a+c_{1,1,b}xyz^b).
\end{eqnarray*}
By using Lemma \ref{eqchar1}, this expression simplifies to
\begin{eqnarray*}
-q^{-3} \bar{\chi}(c_{0,0,a}) \sum_{(y,z) \in (\F_q^\times)^2} \bar{\chi}(z^a).
\end{eqnarray*}
If the order of $\bar{\chi}$ is not a divisor of $a$, then it follows from Lemma \ref{eqchar0} that the contribution of $\sigma$ is equal to $0$.
\\ \\
\textsc{Fact 4}: Let $\sigma:=PQ$ with $P(\cdot,\cdot,0)$ and $Q(\cdot,\cdot,0)$, and let $\tau:=PQR$ with $R(\cdot,\cdot,1)$ be the facet not contained in $\{z=0\}$ that contains $\sigma$, then
$\sigma$ and $\tau$ cancel each other out. Indeed, by Lemma \ref{eqchar1} with $f=f_{\sigma}$ it follows that $L_\sigma=(1-q)L_{\tau}$. As mult$(\Delta_\sigma)=1$, we find that $L_\sigma S(\Delta_\sigma)+ L_\tau S(\Delta_\tau)=0$.
\\ \\
\textsc{Fact 5}: Let $\sigma:=PQ$ with $P(\cdot,\cdot,0)$ and $Q(\cdot,\cdot,1)$, then again by Lemma \ref{eqchar1} one finds $L_P=(1-q)L_{\sigma}$. Now let $\tau_1$ and $\tau_2$ be the facets containing $\sigma$ and let $\tau_0$ be the facet in $\{z=0\}$ containing the vertex $P$. With $\delta_P$ the cone (strictly positively) spanned by $\Delta_{\tau_0},\Delta_{\tau_1}$ and$\Delta_{\tau_2}$ we then find that
$L_\sigma S(\Delta_\sigma)+ L_P S(\delta_P)=0$.
\\ \\
\textsc{Fact 6}: Let $\sigma:=PQ$ with $P(\cdot,\cdot,0)$ and $Q(\cdot,\cdot,1)$, and let $\tau_1$ be a non-compact $B_1$-facet containing $\sigma$. Let $\tau_2$ be the non-compact facet containing the vertex $Q$ and sharing a half line with $\tau_1$. Lemma \ref{eqchar0} implies that $\tau_1 \cap \tau_2$ does not contribute in the formula for
$Z_{f}(\chi,s)$.
\\ \\
\textsc{Fact 7}: Let $\sigma:=PQ$ with $P(\cdot,\cdot,0)$ and $Q(\cdot,\cdot,1)$, and let $\tau_1$ be a non-compact $B_1$-facet containing $\sigma$. Let $\tau_0$ be the non-compact facet containing the vertex $P$ and sharing a half line $\sigma_1$ with $\tau_1$. As $L_{\sigma_1}=(1-q)L_{\tau_1}$ and mult$(\Delta_{\sigma_1})=1$, it follows that the contributions of $\tau_1$ and $\sigma_1$ cancel each other out.
\\ \\
From these facts one can derive the contributions of all possible configurations of $B_1$-facets. We begin with the configuration studied (in the local case over $\Q_p$) in
\cite[Proposition 9.6]{BV15} that we mentioned at the beginning of this section.
\\ \\
\textsc{Case 1: the candidate pole is contributed by a configuration of $B_1$-facets in which no two facets that share a $1$-dimensional face are $B_1$ only for different variables.}
\\ \\
For the contributions to the local Igusa zeta function, one can derive from Facts 1, 4 and 5
that the candidate pole is a false pole. For the global Igusa zeta function, one in addition uses Facts 6 and 7.
\\ \\
\textsc{Case 2: the candidate pole is contributed by exactly two compact $B_1$-facets w.r.t.\ different variables, having a line segment in common.}
\\ \\
If the common line segment is compact, then the configuration is as in Figure 1
%\vspace*{-0.5cm}
\begin{figure}[h]
\begin{center}
\resizebox{2.5in}{!}{\includegraphics{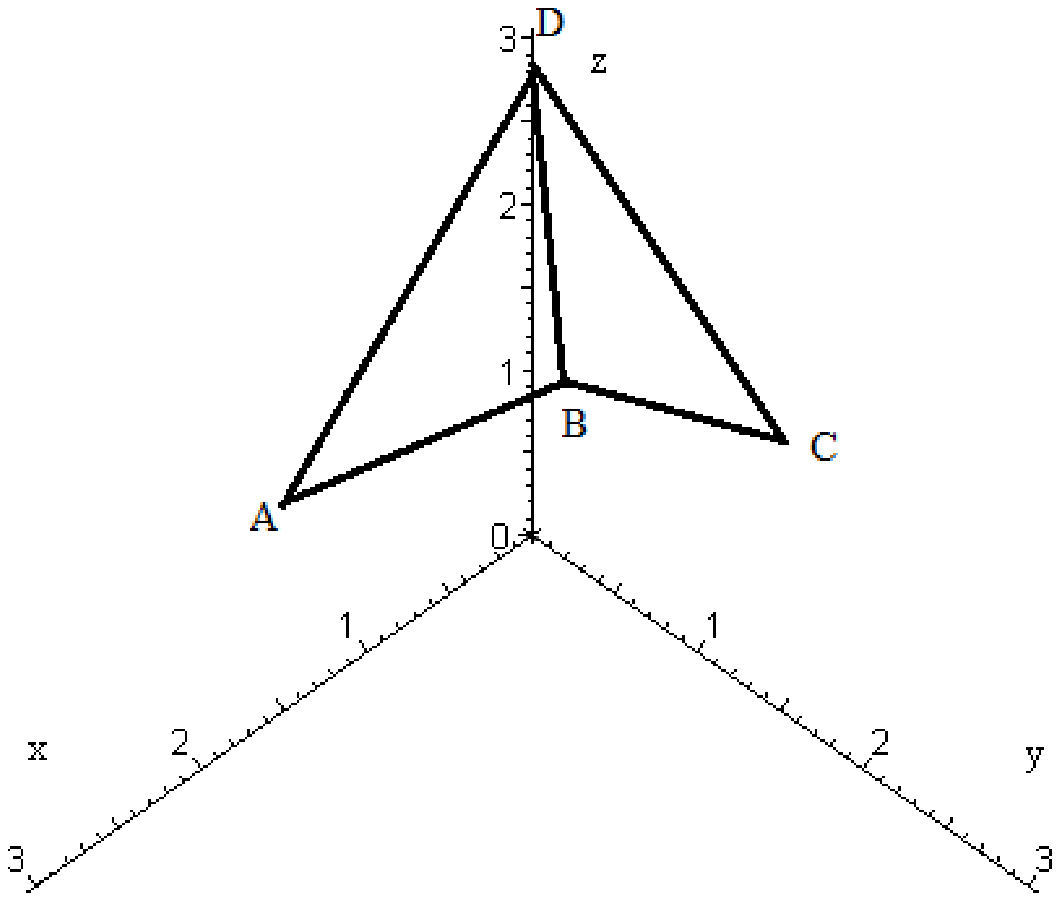}}
\\ Figure 1 \vspace*{-0.5cm}
\end{center}
\end{figure}
with $A(.,0,.), B(1,1,b), C(0,.,.)$ and $D(0,0,a)$.
If the order of $\bar{\chi}$ is not a divisor of $a$, then it follows from Facts 1 to 5 that the candidate pole is a false pole of $Z_{f,0}(\chi,s)$ and $Z_{f}(\chi,s)$.
%Let $\tau_1$ be the compact facet with vertices $A, B$ and $D$ and let $\tau_2$ be the compact facet with vertices $B, C$ and $D$. To compute the contribution of $\sigma=BD$, we consider
%\begin{eqnarray*}
%L_\sigma=p^{-3}\sum_{(x,y,z) \in (\F_p^\times)^3} \bar{\chi}(z^a+xyz^b).
%\end{eqnarray*}
%By using Lemma \ref{eqchar1}, this expression simplifies to
%\begin{eqnarray*}
%-p^{-3}\sum_{(y,z) \in (\F_p^\times)^2} \bar{\chi}(z^a).
%\end{eqnarray*}
%If the order of $\bar{\chi}$ is not a divisor of $a$, then it now follows from Lemma \ref{eqchar0} that the contribution of $\sigma$ is equal to $0$.
%Analogously, if the order of $\bar{\chi}$ is not a divisor of $a$, then the vertex $D$ does not contribute.
%One easily verifies, again by using Lemma \ref{eqchar1}, that the contributions of the segment $AD$ and $\tau_1$ cancel each other.
%As $\chi$ is not the trivial character, one immediately deduces from Lemma \ref{eqchar0} that the contribution of $B$ is equal to $0$.
%To compute the contribution of $\tau=A$, one takes a simplicial subdivision of $\Delta_{\tau}$.
%Let $\tau_3$ be the other facet containing the segment $AB$, denote the facet in $\{y=0\}$ containing $AD$ by $\tau_0$ and let $\delta_A=<\Delta_{\tau_0},\Delta_{\tau_1},\Delta_{\tau_3}>$. Then
%again from Lemma \ref{eqchar1} one deduces that the contributions of $\delta_A$ and $AB$ cancel each other.
%
%In the facet $\tau_2$ one gets the same cancelations which show that the candidate pole is a false pole of $Z_{f,0}(\chi,s)$ and $Z_{f}(\chi,s)$ whenever the order of $\bar{\chi}$ is not a divisor of $a$.
${}$
\\ \\
\textsc{Case 3: the candidate pole is contributed by two non-compact $B_1$-facets w.r.t.\ different variables, having a line segment in common.}
\\ \\
If the common line segment is non-compact, then the configuration is as in Figure 2,
\begin{figure}[h]
\begin{center}
\resizebox{2.5in}{!}{\includegraphics{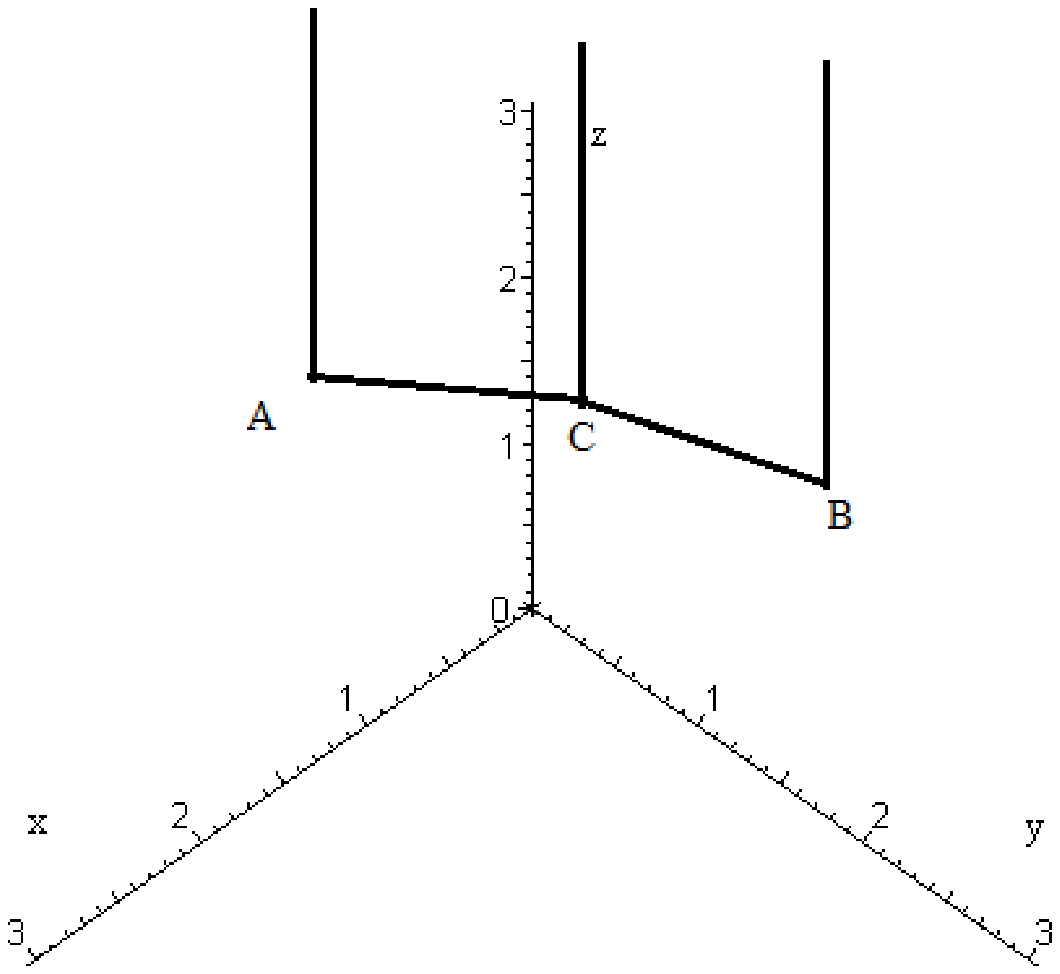}}
\\ Figure 2 \vspace*{-0.5cm}
\end{center}
\end{figure}
with $A(.,0,.), B(0,.,.)$ and $C(1,1,.)$.
For the contributions to the local Igusa zeta function, one deduces from Facts 1, 4 and 5
that the candidate pole is a false pole. For the global Igusa zeta function, one also has to use Facts 6 and 7.

If the common line segment is compact, then its vertices are given by $A(0,0,a)$ and $B(1,1,b)$. If the order of $\bar{\chi}$ is not a divisor of $a$, then by Facts 1 to 7 it follows again that the candidate pole is not an actual pole of $Z_{f,0}(\chi,s)$ and $Z_{f}(\chi,s)$.
${}$
\\ \\
\textsc{Case 4: the candidate pole is contributed by one compact $B_1$-facet and one non-compact $B_1$-facets w.r.t.\ different variables, having a line segment in common.}
\\ \\
Again using Fact 1 to Fact 7, one finds that the candidate pole is a false pole of $Z_{f,0}(\chi,s)$ and $Z_{f}(\chi,s)$ when the order of $\bar{\chi}$ is not a divisor of $a$,.
${}$
\\ \\
\textsc{Case 5: the candidate pole is contributed by at least two $B_1$-facets w.r.t.\ different variables, having a line segment in common.}
\\ \\
As in Case 2 the contributions of $\tau_1:=ABD$, $\tau_2:=BCD$ and $\tau_1 \cap \tau_2:=BD$ are all equal to $0$, one can deduce the fakeness of the candidate pole also when there are other $B_1$-facets having a $1$-dimensional intersection with $\tau_1$ or $\tau_2$. \hfill $\blacksquare$
\\ \\
We now treat the candidate poles contributed by $X_2$-facets.
\begin{lm} \label{lemmaX2}
Let $\tau$ be a facet with vertices $p(N(\tau)/2,0,0), q(x_1,0,2)$ and $r(x_2,2,0)$ where $N(\tau)/2 - x_1$ and $N(\tau)/2 - x_2$ are odd. If the order of $\bar{\chi}$ does not divide $N(\tau)/2$ and is different from $2$, then $\tau$ does not contribute to $Z_{f,0}(\chi,s)$ and $Z_{f}(\chi,s)$.
\end{lm}
\noindent \emph{Proof.}
%Suppose that $-\nu(\tau)/N(\tau)$ is the real part of a candidate pole of $Z_{f,0}(\chi,s)$ or $Z_{f}(\chi,s)$.
%%Now let $d$ be the order of $\chi$ and.
%%If $d | N(\tau)/2$, then $d$ also divides the order of an eigenvalue of monodromy. Indeed, $e^{-\frac{2\pi i}{N(\tau)/2}}$ is a zero of $F_\tau=(1-t^{N(\tau)/2})$.
%If the order of $\bar{\chi}$ does not divide $N(\tau)/2$ and is different from $2$, then we will show that $-\nu(\tau)/N(\tau)$ is not the real part of an actual pole.
It follows immediately from Fact~2 that the vertices $p$, $q$ and $r$ do not contribute. Using
Lemma~\ref{vanishingofLtau} one also verifies that the edge $qr$ does not contribute.
We now show that the contributions of $\sigma_1:=pq$, $\sigma_2:=pr$ and the facet $\tau$ cancel each other.
As $N(\sigma_1)=N(\sigma_2)=N(\tau)$, we have that mult$(\Delta_{\sigma_1})= \text{mult}(\Delta_{\sigma_2}) =1$, and thus
$$S(\Delta_{\sigma_i})=\frac{1}{(q-1)(q^{N(\tau)s+\nu(\tau)} - 1)}, \quad 1 \leq i \leq 2.$$ One gets
\begin{center}
$L_{\sigma_1}S(\Delta_{\sigma_1}) + L_{\sigma_2}S(\Delta_{\sigma_2}) + L_{\tau}S(\Delta_{\tau}) = 0$ \\
$\Updownarrow$ \\
$(q-1)L_\tau = -L_{\sigma_1} - L_{\sigma_2}.$
\end{center}
The equality between these character sums is proven in Proposition~\ref{eqchar2}. \hfill $\blacksquare$\\

If $x_1=x_2=0$ (in which case the $X_2$-facet is the only compact facet of $\Gamma_0$) we can prove something slightly stronger.

\begin{lm} \label{lemmaX2ord2}
Let $\tau$ be a facet with vertices $p(N(\tau)/2,0,0), q(0,0,2)$ and $r(0,2,0)$ where $N(\tau)/2$ is odd. If the order of $\bar{\chi}$
does not divide $N(\tau)/2$,
then $\tau$ does not contribute to $Z_{f,0}(\chi,s)$ and $Z_{f}(\chi,s)$.
\end{lm}
\noindent \emph{Proof.}
The previous proof remains valid, except for the conclusions that $q$, $r$ and $\sigma_3 := qr$ do not contribute, where we used that
the order of $\bar{\chi}$ is not $2$. We show that the contributions cancel.
Indeed, since $\text{mult}(\Delta_q) = \text{mult}(\Delta_r) = \text{mult}(\Delta_{\sigma_3}) = N(\tau)/2$ we have
$$S(\Delta_{\sigma_3})=\frac{N}{(q-1)(q^{N(\tau)s+\nu(\tau)} - 1)}, \quad S(\Delta_q)=S(\Delta_r) = \frac{N}{(q-1)^2(q^{N(\tau)s+\nu(\tau)} - 1)}$$
for some common numerator $N$.
One gets
\begin{center}
$L_qS(\Delta_q) + L_rS(\Delta_r) + L_{\sigma_3}S(\Delta_{\sigma_3}) = 0$ \\
$\Updownarrow$ \\
$(q-1)L_{\sigma_3} = -L_q - L_r.$
\end{center}
This again follows from Proposition~\ref{eqchar2} (with $\alpha = 0$).
\hfill $\blacksquare$

\subsection{The holomorphy conjecture for nondegenerate surface singularities}

\begin{te}
Let $F$ be a number field and let $f(x,y,z) \in \mathcal{O}_F[x,y,z]$ be a polynomial which is nondegenerate over $\C$ w.r.t.\ the compact faces (resp.\ the faces) of its Newton polyhedron at the origin $\Gamma_0$.
Let $K$ be a non-archimedean completion of $F$ with valuation ring $R$ (with maximal ideal $P$) and residue field $\F_q$, and suppose that $\bar{f} := f \bmod P$ is nondegenerate over $\F_q$ w.r.t.\ the compact faces (resp.\ the faces) of $\Gamma_0$. Let $\chi$ be a non-trivial character of $R^\times$ which is trivial on $1 + P$.
Let $\tau$ be a facet of $\Gamma_0$. If $-\nu(\tau)/N(\tau)$ is the real part of a pole of $Z_{f,0}(\chi,s)$ (resp.\ $Z_{f}(\chi,s)$), then the order of $\chi$ divides the order of an eigenvalue of monodromy at some point of $f^{-1}\{ 0 \}$.
\end{te}

\noindent \emph{Proof.}
As before we assume that all compact facets have been subdivided into simplices, without introducing new vertices. We first suppose that $\tau$ is such a simplex. If every $1$-dimensional V-face of $\Gamma_0$ is contained in a compact facet, then we know from Formula (\ref{formulan}) that the zeta function of monodromy at the origin is a product of polynomials. If $\tau$ is not of type $B_1$ nor of type $X_2$, then Proposition \ref{prdeler} implies that the order of $\chi$ divides the order of an eigenvalue of monodromy of $f^{-1}\{ 0 \}$ at the origin.

If $\tau$ is of type $B_1$, then we found in Subsection \ref{subsB1} that there is a point $p(0,0,a)$ in the configuration that is not the intersection of two $1$-dimensional V-faces in a same compact facet, and secondly that the order of $\bar{\chi}$ divides this $a$. This means that the factor $1-t^a$ appears in $\zeta_{f,0}(t)$ and so one finds that the order of $\chi$ divides the order of some eigenvalue of monodromy of $f^{-1}\{ 0 \}$ at the origin.

If $\tau$ is of type $X_2$, say with vertices $p(N(\tau)/2,0,0), q(x_1,0,2)$ and $r(x_2,2,0)$, then $F_\tau=1-t^{N(\tau)/2}$ and hence
$e^{-2\pi i/(N(\tau)/2)}$ is an eigenvalue of monodromy of $f^{-1}\{ 0 \}$ at the origin. Thus if the order of $\bar{\chi}$ divides $N(\tau)/2$ then we are done.
If the order of $\bar{\chi}$ does not divide $N(\tau)/2$, then by Lemma \ref{lemmaX2}
 the order of $\bar{\chi}$ should be equal to $2$. In this situation $N(\tau)/2$ is odd and $x_1$ and $x_2$ are even, while by Lemma~\ref{lemmaX2ord2} we can assume that $0 \neq x_1 \geq x_2$.
 Let $\tau'$ be the other facet which contains the segment $qr$. Notice that $N(\tau')$ is even and that $\tau'$ is not of type $B_1$.
We first suppose that $\tau'$ is compact.
If $\tau'$ is not of type $X_2$, then it follows from Proposition \ref{prdeler} that $e^{-2\pi i/N(\tau')}$ is a zero of $F_{\tau'}$ and so the order of $\chi$ divides the order of some eigenvalue of monodromy of $f^{-1}\{ 0 \}$ at the origin.
If $\tau'$ is of type $X_2$, then the configuration is as in Figure 3.
\begin{figure}[h]
\begin{center}
\resizebox{2.5in}{!}{\includegraphics{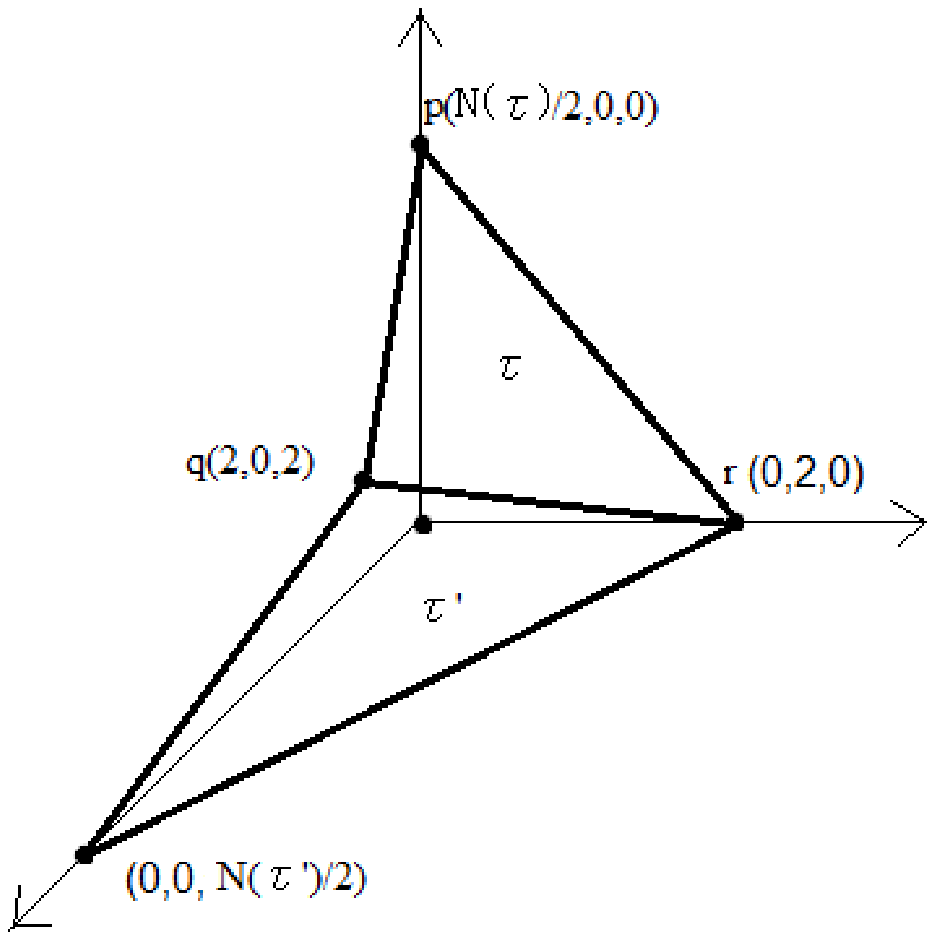}}
\\ Figure 3 \vspace*{-0.5cm}
\end{center}
\end{figure}

In this situation, we get
$$\zeta_{f,0}(t)=(1-t^{N(\tau)/2})(1-t^{N(\tau')/2})(1-t^2),$$
and so again the order of $\chi$ divides the order of an eigenvalue of monodromy of $f^{-1}\{ 0 \}$ at the origin.

Suppose now that $\tau'$ is not compact, then necessarily $x_1 > x_2$ and
$$\mbox{Aff}(\tau') \leftrightarrow x + \frac{x_1-x_2}{2}y = x_1.$$
At a generic point $(0,0,c)$ of the hypersurface, the polynomial $g(x,y,z):=f(x,y,z-c)$ is still nondegenerate w.r.t.\ the compact faces of its Newton polyhedron at the origin (see Lemma \ref{nondegeneratestable}) and its Newton polyhedron is the projection onto $\{z=0\}$ of the Newton polyhedron of $f$ times $\R_+$.
From Varchenko's formula one sees that this projected polyhedron fully determines $\zeta_{g,0}(t)$. Using \cite[Proposition 5]{LVP11} it
follows that $\zeta_{g,0}(t)$ contains the factor $1/(1-t^{x_1})$. We thus have that the order of $\chi$ divides the order of an eigenvalue of monodromy at a point of the hypersurface in the neighbourhood of the origin.

Suppose now that there is a $1$-dimensional V-face $\sigma$, say in the coordinate plane $z=0$, which is not contained in a compact facet. If $e^{-2\pi i/N(\tau)}$ is a zero of $F_{\sigma}$ (we use the notation $F_\sigma$ as if $\sigma$ was a facet of a two-dimensional Newton polyhedron in the plane $z= 0$), then we choose $c \in \mathbb{C}$ close to zero such that $g(x,y,z):=f(x,y,z-c)$ is still nondegenerate w.r.t.\ its Newton polyhedron at the origin (see Lemma \ref{nondegeneratestable}). Then we have $\zeta_{g,0}(t)=\prod\limits_{\sigma \text{ compact facet}}F_{\sigma}$, with $F_{\sigma}=1/\mbox{polynomial}$ (except the case where $\sigma$ contains two vertices on coordinate axes, but in this case the same conclusion holds) and so we find that $e^{-2\pi i/N(\tau)}$ is an eigenvalue of monodromy of $f$ at $(0,0,c)$.
%Suppose now that the configuration is as in Figure 2 and that the order of $\chi$ divides $a$. If $e^{-\frac{2\pi i}{a}}$ is also a zero of $F_\sigma$, then by the same argument we find
%$e^{-\frac{2\pi i}{a}}$ as an eigenvalue of monodromy of $f$ at some point close to the origin.

Finally let $\tau$ be non-compact. Again by the nondegeneracy argument (Lemma \ref{nondegeneratestable}) we can reduce the dimension and conclude that $e^{-2\pi i/N(\tau)}$ is an eigenvalue of monodromy of $f$ at a point in the neighbourhood the origin.
\hfill $\blacksquare$

%\begin{rem} \label{remgeneral}
%\emph{The proof of \cite[Theorem 3.4 ]{H02} can also be written down in the general context where $f$ is a polynomial over an arbitrary number
%field $F$ in $\C$ and where $K$ is a completion w.r.t.\ a finite prime.
%Our proofs in Section \ref{sectioncharsums} and \ref{subsB1} remain valid over the residue field $\F_q$ and hence the holomorphy conjecture for
%nondegenerate surface singularities holds in the general context.}
%\end{rem}

\noindent
\textsc{Vakgroep Wiskunde, Universiteit Gent\\
Krijgslaan 281, 9000 Gent, Belgium}\\

\vspace{-0.4cm}
\noindent \textsc{Departement Elektrotechniek, KU Leuven\\
Kasteelpark Arenberg 10/2452, 3001 Leuven, Belgium}\\

\vspace{-0.4cm}
\noindent \emph{E-mail address:} wouter.castryck@gmail.com
\\ \\
\noindent
\textsc{Faculty of Mathematics and Informatics, Ovidius University\\
BD. Mamaia 124, 900527 Constanta, Romania} \\

\vspace{-0.4cm}
\noindent \emph{E-mail address:} denis@univ-ovidius.ro
\\ \\
\noindent
\textsc{Laboratoire Paul Painlev\'e, Universit\'e Lille 1 \\
Cit\'e Scientifique, 59655 Villeneuve d'Ascq Cedex, France}\\

\vspace{-0.4cm}
\noindent \emph{E-mail address:} lemahieu.ann@gmail.com

\end{document}